\documentclass[12pt,a4paper]{amsart}

\usepackage{amssymb}
\usepackage{amsrefs}
\usepackage[all]{xy}

\allowdisplaybreaks

\newcommand{\bbA}{\mathbb{A}}
\newcommand{\bbK}{\mathbb{K}}
\newcommand{\bbL}{\mathbb{L}}
\newcommand{\bbN}{\mathbb{N}}
\newcommand{\bbZ}{\mathbb{Z}}

\newcommand{\calB}{\mathcal{B}}
\newcommand{\calC}{\mathcal{C}}
\newcommand{\calK}{\mathcal{K}}
\newcommand{\calL}{\mathcal{L}}
\newcommand{\calP}{\mathcal{P}}
\newcommand{\calS}{\mathcal{S}}
\newcommand{\calU}{\mathcal{U}}

\newtheorem*{prop*}{Proposition}
\newtheorem*{theo*}{Theorem}
\newtheorem{lemm}{Lemma}[subsection]

\let\mod=\undefined
\DeclareMathOperator{\Id}{Id} %
\DeclareMathOperator{\Hom}{Hom} %
\DeclareMathOperator{\Ker}{Ker} %
\DeclareMathOperator{\mod}{mod} %

\newcommand{\ol}{\overline}

\title[Normal forms of modules over two-ray algebras]%
  {Normal forms of modules over admissible algebras with formal two-ray modules}

\author{Grzegorz Bobi\'nski}

\address{Faculty of Mathematics and Computer Science \\
Nicolaus Copernicus University \\ Chopina 12/18 \\ 87-100 Toru\'n
\\ Poland}

\email{gregbob@mat.uni.torun.pl}



\date{}

\begin{document}

\begin{abstract}
The aim of the paper is to classify the indecomposable modules and
describe the Auslander--Reiten sequences for admissible algebras
with formal two-ray modules.
\end{abstract}

\maketitle

\section*{Introduction}

Throughout the paper $k$ is a fixed algebraically closed field.
All considered categories are additive $k$-categories and all
functors are $k$-functors.

One of the aims of the representation theory of finite-dimensional
algebras is a description of indecomposable modules and
homomorphism spaces between them. A guiding example is that of
special biserial algebras, for which a full description of the
indecomposable modules and the Auslander--Reiten sequences was
given by Wald and Waschb\"usch~\cite{WalWas} (see
also~\cite{BuRi}). Homomorphism spaces between indecomposable
modules were also investigated (see for example~\cite{CB1}).
Another class of algebras whose representation theory is described
is formed by clannish algebras (or more generally, clan problems)
introduced by Crawley-Boevey~\cite{CB4} (see also~\cites{Bond,
De}). Homomorphism spaces and Auslander--Reiten sequences for this
class of problems were studied by Gei\ss~\cite{Ge} (see
also~\cite{GePe} for a description of the Auslander-Reiten
components).

According to Drozd's Tame and Wild Theorem~\cite{Dr} (see
also~\cite{CB2}) one may hope to obtain classifications like these
above only for so called tame algebras. First examples of tame
algebras are provided by the representation-finite algebras, for
which there are only finitely many isomorphism classes of
indecomposable modules. The representation theory of the
representation-finite algebras has been intensively studied (see for
example~\cites{BaGaRoSa, Bong, BongGa, BrGa}) and seems to be
well-understood. One knows that an algebra is representation-finite
if and only if its infinite radical vanishes.

The first level in the hierarchy of representation-infinite algebras
is occupied by the domestic algebras, for which in each dimension
all but finitely many indecomposable modules can be parameterized by
finitely many lines (see also~\cite{CB3} for a different
characterization of the domestic algebras). Schr\"oer's
work~\cite{Sc} on the infinite radical of special biserial algebras
gives hope to characterize the domestic algebras in terms of the
infinite radial. In~\cite{BobDrSk} (continued by~\cites{Bob1,
BobSk2}), we initiated the study of a new class of domestic
algebras, which may be seen as a test class for this
characterization. The results obtained so far concern the
Auslander--Reiten theory. In order to deal with the infinite radical
one needs to have a more precise knowledge about indecomposable
modules and homomorphisms spaces between them. In this paper we make
a first step in this direction, namely we give a description of the
indecomposable modules. This description resembles the description
obtained for clans, thus one may hope that the corresponding results
about homomorphisms can be also transferred.

The paper is organized as follows. In Section~\ref{mainres} we
present the main result of the paper, in Section~\ref{sectvect} we
recall necessary information about vector space categories, and in
final Section~\ref{sectproof} we prove the main theorem. The paper
was written during the author held a one year post-doc position at
the University of Bern. Author gratefully acknowledges the support
from the Schweizerischer Nationalfonds and the Polish Scientific
Grant KBN No.~1 P03A 018 27.

\section{Strings, the corresponding modules and the main result}
\label{mainres}

In this section we first introduce notation, which is necessary to
formulate the main result of the paper given at the end of the
section.

\subsection{}
In the paper, by $\bbZ$ (respectively, $\bbN_0$, $\bbN$) we denote
the set of (nonnegative, positive) integers. If $m$ and $n$ are
integers, then by $[m, n]$ we denote the set of all integers $l$
such that $m \leq l \leq n$. For a sequence $f : [1, n] \to \bbN$,
$n \in \bbN_0$, of positive integers we denote $n$ by $|f|$. We
identify finite subsets of $\bbN$ with the corresponding
increasing sequences of positive integers. In particular, if $F$
is a finite subset of $\bbN$ and $i \in [1, |F|]$, then $F_i$
denotes the $i$-th element of $F$ with respect to the usual order
of integers.

\subsection{}
By a quiver $Q$ we mean an oriented graph, i.e., a set of vertices
$Q_0$, a set of arrows $Q_1$ and two maps $s_Q, t_Q : Q_1 \to Q_0$,
which assign to an arrow $\alpha$ in $Q$ its starting and
terminating vertex, respectively. If $\alpha \in Q_1$, $s_Q (\alpha)
= x$ and $t_Q (\alpha) = y$, then we write $\alpha : x \to y$. By a
path in $Q$ we mean a sequence $\rho = \alpha_1 \cdots \alpha_n$ of
arrows in $Q$ such that $t_Q (\alpha_{i + 1}) = s_Q (\alpha_i)$ for
all $i \in [1, n - 1]$. The number $n$ is called the length of
$\rho$ and denoted $|\rho|$. We write $s_Q (\rho)$ for $s_Q
(\alpha_n)$ and $t_Q (\rho)$ for $t_Q (\alpha_1)$, and we say that
$\rho$ starts at $s_Q (\rho)$ and terminates at $t_Q (\rho)$. For
each vertex $x$ of $Q$ we denote also by $x$ the path of length $0$
at vertex $x$ ($s_Q (x) = x = t_Q (x)$). For paths $\rho = \alpha_1
\cdots \alpha_n$ and $\rho' = \alpha_1' \cdots \alpha_m'$ in $Q$
such that $s_Q (\rho) = t_Q (\rho')$, we denote by $\rho \rho'$ the
path $\alpha_1 \cdots \alpha_n \alpha_1' \cdots \alpha_m'$. In
particular, $\rho s_Q (\rho) = \rho = t_Q (\rho) \rho$.

\subsection{}
By a defining system we mean a quadruple $(p, q, S, T)$, where $p$
and $q$ are sequences of positive integers such that $|q| = |p|$
and $\sum_{i = 1}^{|p|} p_i \geq 2$, and $S = (S_i)_{i = 1}^{|p|}$
and $T = (T_i)_{i = 1}^{|p|}$ are families of finite subsets of
$\bbN$ such that for each $i \in [1, |p|]$ hold: $T_i \subseteq
S_i \subseteq [2, p_i + |T_i|]$, if $j \in S_i$ then $j + 1 \not
\in S_i$, and $p_i + |T_i| \not \in T_i$. We write $T_{i, j}$
instead of $(T_i)_j$ for $i \in [1, |p|]$ and $j \in [1, |T_i|]$.
Throughout the rest of the section $(p, q, S, T)$ is a fixed
defining system.

\subsection{}
We define a quiver $Q$ by
\begin{align*}
Q_0 & = \{ x_{i, j} \mid i \in [1, |p|], \, j \in [0, p_i + |T_i|]
\} \\
& \cup \{ y_{i, j} \mid i \in [1, |p|], \, j \in [1, q_i - 1] \}
\\
& \cup \{ z_{i, j} \mid i \in [1, |p|], \, j \in S_i \} \\
\intertext{and} %
Q_1 & = \{ \alpha_{i, j} : x_{i, j} \to x_{i, j - 1} \mid i \in
[1, |p|], \, j \in [1, p_i + |T_i|] \} \\
& \cup \{ \beta_{i, j} : y_{i, j} \to y_{i, j - 1} \mid i \in [1,
|p|], \, j \in [1, q_i] \} \\
& \cup \{ \gamma_{i, j} : z_{i, j} \to x_{i, j} \mid i \in [1,
|p|], \, j \in S_i \} \\
& \cup \{ \xi_{i, j} : x_{i, p_i + j} \to z_{i, T_{i, j}} \mid i
\in [1, |p|], \, j \in [1, |T_i|] \},
\end{align*}
where $y_{i, 0} = x_{i + 1, 0}$ (with $x_{|p| + 1, 0} = x_{1, 0}$)
and $y_{i, q_i} = x_{i, p_i}$ for $i \in [1, |p|]$.

Let $A$ be the path algebra of $Q$ bounded by relations
\begin{gather*}
\alpha_{i, j - 1} \alpha_{i, j} \gamma_{i, j}, \, i \in [1, |p|],
\, j \in S_i, \\
\beta_{i, q_i} \alpha_{i, p_i + 1}, \, i \in [1, |p|] \text{ such
that } |T_i| > 0, \\
\xi_{i, j - 1} \alpha_{i, p_i + j}, \, i \in [1, |p|], \, j \in
[2, |T_i|], \\
\intertext{and} %
\alpha_{i, T_{i, j}} \gamma_{i, T_{i, j}} \xi_{i, j} - \alpha_{i,
T_{i, j}} \alpha_{i, T_{i, j} + 1} \cdots \alpha_{i, p_i + j - 1}
\alpha_{i, p_i + j}, \, i \in [1, |p|], \, j \in [1, |T_i|].
\end{gather*}
Recall that by~\cite{Bob1}*{Theorem~1.1} the class of algebras
defined in the above way coincides with the class of admissible
algebras with formal two-ray modules introduced in~\cite{BobSk2}.

In order to clarify a bit the above definitions we give a simple
example. If $p = (6, 3)$, $q = (2, 2)$, $S = (\{ 2, 4, 6, 8 \}, \{ 2
\})$ and $T = (\{ 4, 6 \}, \varnothing)$, then $A$ is the path
algebra of the quiver
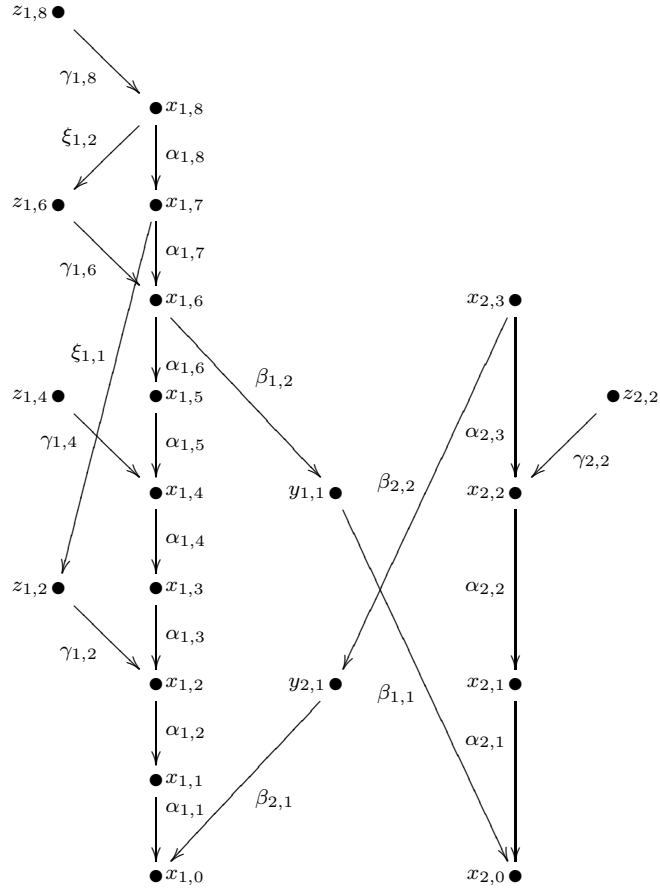
\[
\xymatrix{%
\bullet \save*+!R{\scriptstyle z_{1, 8}} \restore
\ar[rd]_{\gamma_{1, 8}} \\ %
& \bullet \save*+!L{\scriptstyle x_{1, 8}} \restore
\ar[d]^{\alpha_{1, 8}} \ar[ld]_{\xi_{1, 2}} \\ %
\bullet \save*+!R{\scriptstyle z_{1, 6}} \restore
\ar[rd]_{\gamma_{1, 6}} & \bullet \save*+!L{\scriptstyle x_{1, 7}}
\restore \ar[d]^{\alpha_{1, 7}} \ar[ldddd]_(.4){\xi_{1, 1}} \\ %
& \bullet \save*+!L{\scriptstyle x_{1, 6}} \restore
\ar[d]^(.7){\alpha_{1, 6}} \ar[rrdd]^{\beta_{1, 2}} & & & & \bullet
\save*+!R{\scriptstyle x_{2, 3}} \restore
\ar[dd]_(.7){\alpha_{2, 3}} \ar[lldddd]_{\beta_{2, 2}} \\ %
\bullet \save*+!R{\scriptstyle z_{1, 4}} \restore
\ar[rd]_(.3){\gamma_{1, 4}} & \bullet \save*+!L{\scriptstyle x_{1,
5}} \restore \ar[d]^{\alpha_{1, 5}} & & & & & \bullet
\save*+!L{\scriptstyle z_{2, 2}} \restore \ar[ld]^{\gamma_{2, 2}} \\ %
& \bullet \save*+!L{\scriptstyle x_{1, 4}} \restore
\ar[d]^{\alpha_{1, 4}} & & \bullet \save*+!R{\scriptstyle y_{1, 1}}
\restore \ar[rrdddd]_{\beta_{1, 1}} & & \bullet
\save*+!R{\scriptstyle x_{2, 2}} \restore \ar[dd]_{\alpha_{2, 2}} \\ %
\bullet \save*+!R{\scriptstyle z_{1, 2}} \restore
\ar[rd]_{\gamma_{1, 2}} & \bullet \save*+!L{\scriptstyle x_{1, 3}}
\restore \ar[d]^{\alpha_{1, 3}} \\ %
& \bullet \save*+!L{\scriptstyle x_{1, 2}} \restore
\ar[d]^{\alpha_{1, 2}} & & \bullet \save*+!R{\scriptstyle y_{2, 1}}
\restore \ar[lldd]^{\beta_{2, 1}} & & \bullet \save*+!R{\scriptstyle
x_{2, 1}} \restore \ar[dd]_(.3){\alpha_{2,
1}} \\ %
& \bullet \save*+!L{\scriptstyle x_{1, 1}} \restore
\ar[d]^(.3){\alpha_{1, 1}} \\ %
& \bullet \save*+!L{\scriptstyle x_{1, 0}} \restore & & & & \bullet
\save*+!R{\scriptstyle x_{2, 0}} \restore}
\]
bounded by relations
\begin{gather*}
\alpha_{1, 1} \alpha_{1, 2} \gamma_{1, 2}, \; \alpha_{1, 3}
\alpha_{1, 4} \gamma_{1, 4}, \; \alpha_{1, 5} \alpha_{1, 6}
\gamma_{1, 6}, \; \alpha_{1, 7} \alpha_{1, 8} \gamma_{1, 8}, \;
\alpha_{2, 1} \alpha_{2, 2} \gamma_{2, 2}, \; \beta_{1, 2}
\alpha_{1, 7},
\\ %
\xi_{1, 1} \alpha_{1, 8}, \; \alpha_{1, 2} \alpha_{1, 3} \alpha_{1,
4} \alpha_{1, 5} \alpha_{1, 6} \alpha_{1, 7} - \alpha_{1, 2}
\gamma_{1, 2} \xi_{1, 1}, \; \alpha_{1, 6} \alpha_{1, 7} \alpha_{1,
8} - \alpha_{1, 6} \gamma_{1, 6} \xi_{1, 2}.
\end{gather*}

\subsection{}
Let
\[
Q_1' = \{ \alpha_{i, j} : x_{i, j} \to x_{i, j - 1} \mid i \in [1,
|p|], \, j \in [1, p_i + |T_i|] \}
\]
and $Q_1'' = Q_1 \setminus Q_1'$. Let $Q^*$ be the quiver with
same set of vertices and arrows as $Q$, but with the arrows from
$Q_1''$ reversed, i.e., $Q_0^* = Q_0$, $Q_1^* = Q_1$ and
\[
s_{Q^*} (\alpha) =
\begin{cases}
s_Q (\alpha) & \alpha \in Q_1', \\
t_Q (\alpha) & \alpha \in Q_1'',
\end{cases} \text{ and }
t_{Q^*} (\alpha) =
\begin{cases}
t_Q (\alpha) & \alpha \in Q_1', \\
s_Q (\alpha) & \alpha \in Q_1''.
\end{cases}
\]

By a string in $Q$ we mean a path in $Q^*$ which does not contain
a subpath $\alpha_{i, T_{i, j}} \alpha_{i, T_{i, j} + 1} \cdots
\alpha_{i, p_i + j}$ for $i \in [1, |p|]$ and $j \in [1, |T_i|]$.
For formal reasons we also introduce the empty string denoted by
$\varnothing$. By convention the length of $\varnothing$ is $-1$,
the maps $s_{Q^*}$ and $t_{Q^*}$ are not defined for $\varnothing$
and it cannot be composed with other strings. If $C$ is a string
and $C = C' C''$ for strings $C'$ and $C''$, then $C'$ is called a
terminating substring of $C$ and $C''$ is called a starting
substring of $C$.

If $C = c_1 \cdots c_n$ is a string and $x \in Q_0$, then we put
\begin{align*}
J_C^x & = \{ i \in [0, n - 1] \mid t_{Q^*} (c_{i + 1}) = x \} \\
\intertext{and} %
I_C^x & =
\begin{cases}
J_C^x \cup \{ n \} & s_{Q^*} (c_n) = x, \\
J_C^x & s_{Q^*} (c_n) \neq x.
\end{cases}
\end{align*}
In particular, $J_y^x = \varnothing$ for all $y \in Q_0$, $I_x^x =
\{ 0 \}$, and $I_y^x = \varnothing$ if $y \neq x$.

\subsection{}
For each vertex $x$ of $Q$ we denote by $\omega_x$ (respectively
$\mu_x$) the longest string terminating at $x$ and consisting only
of elements of $Q_1'$ ($Q_1''$). Similarly, by $\pi_x$
(respectively $\nu_x$) we denote the longest string starting at
$x$ and consisting only of elements of $Q_1'$ ($Q_1'')$.

Let
\begin{align*}
Q_0' = \{ x_{i, j} \mid i \in [1, |p|], \, j \in S_i \}, \\
\intertext{and} %
Q_0'' = \{ x_{i, j} \mid i \in [1, |p|], \, j \in T_i \}.
\end{align*}
For $x \in Q_0'$, $x = x_{i, j}$, we denote $\alpha_{i, j}$ by
$\alpha_x$ and $\gamma_{i, j}$ by $\gamma_x$.

Let $x \in Q_0''$, $x = x_{i, T_{i, j}}$. We put
\[
B_x = \alpha_{i, T_{i, j} + 1} \cdots \alpha_{i, j} \xi_{i, j}
\gamma_{i, T_{i, j}}.
\]
For a string $C$ terminating at $x$ we denote by $p_C$ the maximal
integer $p \geq 0$ such that $B_x^p$ is a terminating substring of
$C$, where $B_x^p$ denotes the $p$-fold composition of $B_x$ with
itself (with the convention that $B_x^0 = x$). If $x \in Q_0'
\setminus Q_0''$ then we set $B_x = x$ and $p_C = 0$ for each
string $C$ terminating at $x$.

\subsection{} \label{sectord}
For a given vertex $x$ of $Q$ we introduce a linear order in the
set of all strings terminating at $x$. Let $C$ and $C'$ be two
strings terminating at $x$ and let $C_0$ be the longest string
which is both a terminating substring of $C$ and a terminating
substring of $C'$. Then $C < C'$ if and only if either $C = C_0
\beta D$ for $\beta \in Q_1''$ and a string $D$ or $C' = C_0
\alpha D'$ for $\alpha \in Q_1'$ and a string $D'$. Note that the
maximal string terminating at $x$ is $\omega_x$ and the minimal
one is $\mu_x$.

If $C \neq \omega_x$ is a string terminating at $x$, then there
exists a direct successor $C_+$ of $C$, which can be described in
the following way. If there exists $\alpha \in Q_1'$ such that $C
\alpha$ is a string, then $C_+ = C \alpha \mu_{s_Q (\alpha)}$.
Otherwise, there exist a string $C'$ and $\beta \in Q_1''$ such
that $C = C' \beta \omega_{t_Q (\beta)}$. In this case $C_+ = C'$.
We also put $(\omega_x)_+ = \varnothing$.

Similarly, we may define a string ${}_+ C$, which is a direct
successor of $C$ with respect to the appropriate order in the set
of all strings starting at $s_{Q^*} (C)$. Since this order will
play no role in the sequel, we only give a description of ${}_+
C$. If there exists $\beta \in Q_1''$ such that $\beta C$ is a
string, then ${}_+ C = \pi_{s_Q (\beta)} \beta C$. Otherwise,
${}_+ C = C''$, if $C = \nu_{t_Q (\alpha)} \alpha C''$ for $\alpha
\in Q_1'$ and a string $C''$, or ${}_+ C = \varnothing$ if $C =
\nu_x$.

Let $C$ be a string such that $|C_+| + |{}_+ C| \geq |C|$ (this is
equivalent to saying that $C \neq \nu_x \omega_x$ for a vertex $x$
of $Q$). Then we define ${}_+ C_+$ by
\[
{}_+ C_+ =
\begin{cases}
{}_+ (C_+) & C_+ \neq \varnothing, \\
({}_+ C)_+ & {}_+ C \neq \varnothing.
\end{cases}
\]
One easily verifies that the above definition is correct and ${}_+
C_+ \neq \varnothing$. We also put ${}_+ (\nu_x \omega_x)_+ =
\varnothing$ for $x \in Q_0$.

\subsection{}
Let $\calS$ be the set of all strings in $Q$. For $x \in Q_0'$ we
denote by $\calS_x$ the set of all strings $C$ terminating at $x$
such that $\alpha_x C'$ is a string, where $C = B_x^{p_C} C'$
($\calS_x$ is the set of all strings terminating at $x$ if $x \in
Q_0' \setminus Q_0''$). Let $\calP_x$ be the set all pairs $(C,
C')$ of $C, C' \in \calS_x$ such that $C < C'$ and, if $x \in
Q_0''$, $C' < B_x C$. Finally, we put
\[
\calB' = \{ B_x \mid x \in Q_0'' \} \text{ and } \calB = \{ B_0 \}
\cup \calB',
\]
where
\[
B_0 = \alpha_{1, 1} \cdots \alpha_{1, p_1} \beta_{1, q_1} \cdots
\beta_{1, 1} \cdots \alpha_{|p|, 1} \cdots \alpha_{|p|, p_{|p|}}
\beta_{|p|, q_{|p|}} \cdots \beta_{|p|, 1}.
\]

\subsection{}
Let $B = b_1 \cdots b_n \in \calB$, $\lambda \in k^*$ and $m \in
\bbN$. We define a representation $R (B, \lambda, m)$ of $Q$ as
follows:
\begin{align*}
R (B, \lambda, m)_y & = \bigoplus_{j \in [1, m]} \bigoplus_{i \in
J_B^y} k v_i^{(j)} \\
\intertext{and} %
R (B, \lambda, m)_\alpha (v_i^{(j)}) & =
\begin{cases}
v_{i - 1}^{(j)} & \alpha \in Q_1', \, \alpha = b_i, \, i \in [1, n
- 1],
\\ %
v_{i + 1}^{(j)} & \alpha \in Q_1'', \, \alpha = b_{i + 1}, \, i
\in [0, n - 2],
\\ %
\lambda v_0^{(j)} + v_0^{(j + 1)} & \alpha = b_n, \, i = n - 1,
\\ %
& \qquad j \in [1, m - 1],
\\ %
\lambda v_0^{(m)} & \alpha = b_n, \, i = n - 1, \, j = m, \\
0 & \text{otherwise}.
\end{cases}
\end{align*}
We also put $R (B, \lambda, 0) = 0$.

\subsection{}
Let $x \in Q_0''$, $B = b_1 \cdots b_n = B_x$ and $m \in \bbN$. We
define a representation $Q (B, m)$ of $Q$ as follows:
\begin{align*}
Q (B, m)_y & =
\begin{cases}
k v' \oplus \bigoplus_{j \in [1, m]} \bigoplus_{i \in J_C^y} k
v_i^{(j)} & y = t_Q (\alpha_x), \\
\bigoplus_{j \in [1, m]} \bigoplus_{i \in J_C^y} k v_i^{(j)} &
\text{otherwise},
\end{cases} \\
Q (B, m)_\alpha (v_i^{(j)}) & =
\begin{cases}
v' & \alpha = \alpha_x, \, i = 0, \, j = 1,
\\ %
v_{i - 1}^{(j)} & \alpha \in Q_1', \, \alpha = b_i, \, i \in [1, n
- 1],
\\ %
v_{i + 1}^{(j)} & \alpha \in Q_1'', \, \alpha = b_{i + 1}, \, i
\in [0, n - 2],
\\ %
v_0^{(j)} + v_0^{(j + 1)} & \alpha = b_n, \, i = n - 1, \, j \in [1, m - 1],
\\ %
v_0^{(m)} & \alpha = b_n, \, i = n - 1, \, j = m,
\\ %
0 & \text{otherwise},
\end{cases}
\intertext{and} %
Q (B, m)_\alpha (v') & = 0.
\end{align*}

\subsection{}
Let $C = c_1 \cdots c_n \in \calS$. We define a representation $M
(C)$ of $Q$ as follows:
\begin{align*}
M (C)_y & = \bigoplus_{i \in I_C^y} k v_i
\\ %
\intertext{and} %
M (C)_\alpha (v_i) & =
\begin{cases}
v_{i - 1} & \alpha \in Q_1', \, \alpha = c_i, \, i \in [1, n],
\\ %
v_{i + 1} & \alpha \in Q_1'', \, \alpha = c_{i + 1}, \, i \in [0,
n - 1],
\\ %
0 & \text{otherwise}.
\end{cases}
\end{align*}
In particular, $M (x)$ is the simple representation of $Q$ at $x$.
We also put $M (\varnothing) = 0$.

\subsection{}
Let $x \in Q_0'$ and $C = c_1 \cdots c_n \in \calS_x$. We define a
representation $N (C)$ of $Q$ as follows:
\begin{align*}
N (C)_y & =
\begin{cases}
k v' \oplus \bigoplus_{i \in I_C^y} k v_i  & y = t_Q (\alpha_x),
\\ %
k v'' \oplus \bigoplus_{i \in I_C^y} k v_i  & y = s_Q (\gamma_x),
\\ %
\bigoplus_{i \in I_C^y} k v_i & \text{otherwise},
\end{cases}
\\ %
N (C)_\alpha (v_i) & =
\begin{cases}
v' & \alpha = \alpha_x, \, i = p |B_x|, \, p \in [0, p_C],
\\ %
v_{i - 1} & \alpha \in Q_1', \, \alpha = c_i, \, i \in [1, n],
\\ %
v_{i + 1} & \alpha \in Q_1'', \, \alpha = c_{i + 1}, \, i \in [0,
n - 1],
\\ %
0 & \text{otherwise},
\end{cases}
\\ %
N (C)_\alpha (v') & = 0,
\\ %
\intertext{and} %
N (C)_\alpha (v'') & =
\begin{cases}
v_0 & \alpha = \gamma_x,
\\ %
0 & \text{otherwise}.
\end{cases}
\end{align*}
We also put $N (\varnothing) = M (s_Q (\gamma_x))$ (more
precisely, we should write $N_x (\varnothing)$,  but we omit the
vertex if it causes no confusion).

\subsection{}
Let $x \in Q_0''$ and $C = c_1 \cdots c_n \in \calS_x$ be such
that $p_C > 0$. We define a representation $L (C)$ of $Q$ as
follows:
\begin{align*}
L (C)_y & =
\begin{cases}
k v' \oplus \bigoplus_{i \in I_C^y} k v_i  & y = t_Q (\alpha_x),
\\ %
\bigoplus_{i \in I_C^y} k v_i & \text{otherwise},
\end{cases}
\\ %
L (C)_\alpha (v_i) & =
\begin{cases}
v' & \alpha = \alpha_x, \, i = p |B_x|, \, p \in [0, p_C],
\\ %
v_{i - 1} & \alpha \in Q_1', \, \alpha = c_i, \, i \in [1, n],
\\ %
v_{i + 1} & \alpha \in Q_1'', \, \alpha = c_i, \, i \in [0, n -
1],
\\ %
0 & \text{otherwise},
\end{cases}
\\ %
\intertext{and} %
L (C)_\alpha (v') & = 0.
\end{align*}

\subsection{}
Let $x \in Q_0'$, and $(C = c_1 \cdots c_n, C' = c_1' \cdots c_m')
\in \calP_x$. We define a representation $N (C, C')$ of $Q$ as
follows:
\begin{align*}
N (C, C')_y & =
\begin{cases}
k v' \oplus \bigoplus_{i \in I_C^y} k v_i \oplus \bigoplus_{i \in
I_{C'}^y} k v_i' & y = t_Q (\alpha_x),
\\ %
k v'' \oplus \bigoplus_{i \in I_C^y} k v_i \oplus \bigoplus_{i \in
I_{C'}^y} k v_i' & y = s_Q (\gamma_x),
\\ %
\bigoplus_{i \in I_C^y} k v_i \oplus \bigoplus_{i \in I_{C'}^y} k
v_i' & \text{otherwise},
\end{cases}
\\ %
N (C, C')_\alpha (v_i) & =
\begin{cases}
v' & \alpha = \alpha_x, \, i = p |B_x|, \, p \in [0, p_C],
\\ %
v_{i - 1} & \alpha \in Q_1', \, \alpha = c_i, \, i \in [1, n],
\\ %
v_{i + 1} & \alpha \in Q_1'', \, \alpha = c_{i + 1}, \, i \in [0,
n - 1],
\\ %
0 & \text{otherwise},
\end{cases}
\\ %
N (C, C')_\alpha (v_i') & =
\begin{cases}
v' & \alpha = \alpha_x, \, i = p |B_x|, \, p \in [0, p_{C'}],
\\ %
v_{i - 1}' & \alpha \in Q_1', \, \alpha = c_i', \, i \in [1, m],
\\ %
v_{i + 1}' & \alpha \in Q_1'', \, \alpha = c_{i + 1}', \, i \in
[0, m - 1],
\\ %
0 & \text{otherwise},
\end{cases} \\ %
N (C, C')_\alpha (v') & = 0,
\\ %
\intertext{and} %
N (C, C')_\alpha (v'') & =
\begin{cases}
v_0 & \alpha = \gamma_x, \\ %
0 & \text{otherwise}.
\end{cases}
\end{align*}
We also put $N (C, \varnothing) = M (\gamma_x C)$,  $N (C, C) = N
(C) \oplus M (C)$ and, if $x \in Q_1''$, $N (C, B_x C) = L (B_x C)
\oplus M (\gamma_x C)$.

\subsection{} \label{maintheo}
Let
\begin{multline*}
\calS' = \calS \setminus (\{ \nu_x \omega_x \mid x \in Q_0 \} \cup
\{ C \mid C \in \calS_x, \, x \in Q_0' \} \\ %
\cup \{ \alpha_x C \mid C \in \calS_x, \, x \in Q_0' \} \cup \{
\gamma_x C \mid C \in \calS_x, \, x \in Q_0'' \}).
\end{multline*}
Observe, that $\nu_x \omega_x \in \calS$ for all $x \in Q_0$, and
$\alpha_x C \in \calS$ for all $x \in Q_0' \setminus Q_0''$ and $C
\in \calS_x$. Moreover, if $x \in Q_0''$ and $C \in \calS_x$, then
$\gamma_x C \in \calS$, but $\alpha_x C \in \calS$ if and only if
$\omega_x$ is not a terminating substring of $C$.

The following theorem is the main result of the paper.

\begin{theo*}
Let $(p, q, S, T)$ be a defining system and let $A$ be the
corresponding algebra.
\begin{enumerate}

\item
Representations
\begin{align*}
& R (B, \lambda, m), \, B \in \calB, \, \lambda \in k^*, \, m \in
\bbN, \\ %
& Q (B, m), \, B \in \calB', \, m \in \bbN, \\ %
& M (C), \, C \in \calS, \\ %
& N (C), \, C \in \calS_x, \, x \in Q_0', \\ %
& L (B_x C), \, C \in \calS_x, \, x \in Q_0'', \\ %
& N (C, C'), (C, C') \in \calP_x, \, x \in Q_0',
\end{align*}
form a complete set of pairwise nonisomorphic indecomposable
modules over $A$.

\item
Sequences
\begin{align*}
& 0 \to R (B, \lambda, m) \to R (B, \lambda, m + 1) \oplus R (B,
\lambda, m - 1) \to R (B, \lambda, m)
\\ %
& \qquad \to 0, \, (B, \lambda, m) \in \calB \times k^* \times
\bbN, \, B = B_0 \text{ or } \lambda \neq 1,
\\ %
& 0 \to R (B, 1, m) \to Q (B, m + 1) \oplus R (B, 1, m - 1) \to Q
(B, m) \to 0,
\\ %
& \qquad B \in \calB', \, m \in \bbN,
\\ %
& 0 \to Q (B, m) \to R (B, m) \oplus Q (B, m - 1) \to R (B, 1, m -
1) \to 0,
\\ %
& \qquad B \in \calB', \, m \in \bbN, \, m > 1,
\\ %
& 0 \to M (C) \to M (C_+) \oplus M ({}_+ C) \to M ({}_+ C_+) \to
0, \, C \in \calS',
\\ %
& 0 \to M (C) \to M (C_+) \oplus N (\mu_x, {}_+ C) \to N (\mu_x,
{}_+ C_+) \to 0,
\\ %
& \qquad C = \alpha_x C', \, C' \in \calS_x, \, x \in Q_0',
\\ %
& 0 \to M (C) \to N (C, C_+) \to N (C_+) \to 0, \, C \in \calS_x,
\, x \in Q_0',
\\ %
& 0 \to M (\gamma_x C) \to N (C_+, B_x C) \to L (B_x C_+) \to 0,
\, C \in \calS_x, \, x \in Q_0'',
\\ %
& 0 \to N (C) \to N (C, C_+) \to M (C_+) \to 0, \, C \in \calS_x,
\, x \in Q_0', \, C \neq \omega_x,
\\ %
& 0 \to L (B_x C) \to N (C_+, B_x C) \to M (\gamma_x C_+) \to 0,
\, C \in \calS_x, \, x \in Q_0'',
\\ %
& 0 \to N (C, C') \to N (C, C_+') \oplus N (C_+, C') \to N (C_+,
C_+') \to 0,
\\ %
& \qquad (C, C') \in \calP_x, \, x \in Q_0',
\end{align*}
form a complete list of Auslander--Reiten sequences in $\mod A$.

\end{enumerate}
\end{theo*}

We finish this section with some remarks concerning the above
theorem. First of all, if $x \in Q_0'$ then $\omega_x \in \calS_x$
if and only if $x \not \in Q_0''$. If $x \in Q_0'$, $C \in
\calS_x$ and $\alpha_x C \in \calS$, then ${}_+ (\alpha_x C) = C$
and ${}_+ (\alpha_x C)_+ = C_+$. Moreover, if $C \neq \omega_x$,
then $(\alpha_x C)_+ = \alpha_x C_+$. Finally, if $x \in Q_0'
\setminus Q_0''$, then $(\alpha_x \omega)_+ = \varnothing$.

\section{Vector space categories} \label{sectvect}

In this section we describe vector space categories and subspace
categories needed in the proof of our main result.

\subsection{}
Following~\cite{Si}*{Section~17.1} (see
also~\cite{Ri2}*{Section~2.4}) by a vector space category we mean
a pair $\bbK = (\calK, {|-|})$, where $\calK$ is a Krull--Schmidt
category and ${|-|} : \calK \to \mod k$ is a faithful functor. For
a vector space category $\bbK$ we consider the subspace category
$\calU (\bbK)$ of $\bbK$. The objects of $\calU (\bbK)$ are
triples $V = (V_0, V_1, \gamma_V)$ with $V_0 \in \calK$, $V_1 \in
\mod k$ and $\gamma_V : V_1 \to |V_0|$ a $k$-linear map. If $V =
(V_0, V_1, \gamma_V)$ and $W = (W_0, W_1, \gamma_W)$ are two
objects of $\calU (\bbK)$, then a morphism $f : V \to W$ in $\calU
(\bbK)$ is a pair $f = (f_0, f_1)$, where $f_0 : V_0 \to W_0$ is a
morphism in $\calK$, $f_1 : V_1 \to W_1$ is a $k$-linear map and
the condition $|f_0| \gamma_V = \gamma_W f_1$ is satisfied. By
$\ol{0}$ we denote the triple $(0, k, 0)$ in $\calU (\bbK)$.

\subsection{}
An ordered set $I$ is called semi-admissible, if the order is linear
and for each element of $I$ which is not maximal there exists a
direct successor. If in addition, there exist a minimal and a
maximal elements in $I$, then we call $I$ admissible. If $I$ is a
semi-admissible ordered set and $\gamma \in I$ is not maximal in
$I$, then by $\gamma_+$ we denote the direct successor of $\gamma$
in $I$.

If $I_1$ and $I_2$ are two semi-admissible ordered sets, then we
introduce the order in $I_1 \times I_2$ by saying that $(x_1, y_1)
\leq (x_2, y_2)$ if either $x_1 < x_2$ or $x_1 = x_2$ and $y_1 \leq
y_2$, for $x_1, x_2 \in I_1$ and $y_1, y_2 \in I_2$. If $(x, y) \in
I_1 \times I_2$, then we put $(x, y)^+ = (x_+, y)$. If in addition
$I_1$ and $I_2$ are disjoint, then by $I_1 + I_2$ we denote the
ordered set $I_1 \cup I_2$ with the elements of $I_1$ smaller than
the elements of $I_2$.

If $I$ is an admissible ordered set, then we denote by $I_-$ the
set $\{ * \} + I$, where $* \not \in I$. Note that in this case $*
= \min I_-$ and $*_+ = \min I$. Similarly, we put $I_+ = I + \{ *
\}$ (thus in this case $* = \max I_+ = (\max I)_+$). Finally, we
denote by $I'$ the ordered set $I \setminus \{ \max I \}$.

\subsection{}
Let $I_1$, \ldots, $I_{r + 1}$, $r \in \bbN_0$, be a family of
admissible ordered sets. Let $\calK$ be the Krull--Schmidt
category, whose indecomposable objects are
\begin{itemize}

\item
$X_\gamma$, $\gamma \in I_p'$, $p \in [1, r + 1]$,

\item
$X_{\max I_p}'$, $X_{\max I_p}''$, $p \in [1, r]$,

\end{itemize}
and all indecomposable objects of $\calK$ are one-dimensional, i.e.,
for each indecomposable object $X$ of $\calK$,  $\dim_k |X| = 1$. If
$U$ and $V$ are indecomposable objects of $\calK$, then
$\Hom_{\calK} (U, V) \neq 0$ if and only if one of the following
conditions holds:
\begin{itemize}

\item
$U = X_{\gamma'}$, $V = X_{\gamma''}$, $\gamma' \in I_p'$,
$\gamma'' \in I_q'$, $(p, \gamma') \leq (q, \gamma'')$,

\item
$U = X_\gamma$, $V = X_{\max I_q}'$, $\gamma \in I_p'$, $p \leq
q$,

\item
$U = X_\gamma$, $V = X_{\max I_q}''$, $\gamma \in I_p'$, $p \leq
q$,

\item
$U = X_{\max I_p}'$, $V = X_\gamma$, $\gamma \in I_q'$, $p < q$,

\item
$U = X_{\max I_p}'$, $V = X_{\max I_q}'$, $p \leq q$,

\item
$U = X_{\max I_p}'$, $V = X_{\max I_q}''$, $p < q$,

\item
$U = X_{\max I_p}''$, $V = X_\gamma$, $\gamma \in I_q'$, $p < q$,

\item
$U = X_{\max I_p}''$, $V = X_{\max I_q}'$, $p < q$,

\item
$U = X_{\max I_p}''$, $V = X_{\max I_q}''$, $p \leq q$.

\end{itemize}
By $\bbK_{I_1, \ldots, I_{r + 1}}$ we denote the vector space
category $(\calK, {|-|})$, where ${|-|} : \calK \to \mod k$ is the
forgetful functor.

\subsection{}
Let $I$ be an admissible ordered set. Let $\calL$ be the
Krull--Schmidt category, whose indecomposable objects are
\begin{itemize}

\item
$X_\gamma$, $\gamma \in I$,

\item
$Y_\gamma$, $\gamma \in I$,

\end{itemize}
and all indecomposable objects of $\calL$ are one-dimensional. If
$U$ and $V$ are indecomposable objects of $\calL$, then
$\Hom_{\calL} (U, V) \neq 0$ if and only if one of the following
conditions holds:
\begin{itemize}

\item
$U = X_{\gamma'}$, $V = X_{\gamma''}$, $\gamma' \leq \gamma''$,

\item
$U = X_{\gamma'}$, $V = Y_{\gamma''}$, $\gamma' \leq \gamma''$,

\item
$U = Y_{\gamma'}$, $V = Y_{\gamma''}$, $\gamma' \leq \gamma''$.

\end{itemize}
By $\bbL_I$ we denote the vector space category $(\calL, {|-|})$,
where ${|-|} : \calL \to \mod k$ is the forgetful functor.

\subsection{} \label{subspaceone}
We have the following description of the indecomposable objects and
the Auslander--Reiten sequences in $\calU (\bbL_I)$. For definitions
of the relevant objects and the proof we refer
to~\cite{BobDrSk}*{Section~3}.

\begin{prop*}
Let $I$ be an admissible ordered set.
\begin{enumerate}

\item
Objects
\begin{align*}
& M_{\min I_-, \gamma} = X_\gamma, \, \gamma \in I, \,
\\ %
& M_{\gamma', \gamma''} = \ol{Y_{\gamma'} X_{\gamma''}}, \,
\gamma', \gamma'' \in I, \, \gamma' < \gamma'',
\\ %
& M_{\gamma, \max I_+} = \ol{Y_\gamma}, \, \gamma \in I,
\\ %
& M_{\gamma, \gamma}' = Y_\gamma, \, \gamma \in I,
\\ %
& M_{\gamma, \gamma}'' = \ol{X_\gamma}, \, \gamma \in I,
\\ %
& M_{\max I_+, \max I_+}'' = \ol{0},
\end{align*}
form a complete set of pairwise nonisomorphic indecomposable
objects in $\calU (\bbL_I)$.

\item
Sequences
\begin{align*}
& 0 \to M_{\gamma', \gamma''} \to M_{\gamma'_+, \gamma''} \oplus
M_{\gamma', \gamma''_+} \to M_{\gamma'_+, \gamma''_+} \to 0, \,
\gamma', \gamma'' \in I_-, \, \gamma' < \gamma'',
\\ %
& 0 \to M_{\gamma, \gamma}' \to M_{\gamma, \gamma_+} \to
M_{\gamma_+, \gamma_+}'' \to 0, \, \gamma \in I,
\\ %
& 0 \to M_{\gamma, \gamma}'' \to M_{\gamma, \gamma_+} \to
M_{\gamma_+, \gamma_+}' \to 0, \, \gamma \in I',
\end{align*}
form a complete list of Auslander--Reiten sequences in $\calU
(\bbL_I)$, where
\begin{align*}
& M_{\gamma, \gamma} = M_{\gamma, \gamma}' \oplus M_{\gamma,
\gamma}'', \, \gamma \in I,
\\ %
& M_{\min I_-, \max I_+} = 0.
\end{align*}

\end{enumerate}
\end{prop*}

\subsection{}
Let $I_0$, \ldots, $I_{r + 1}$, $r \in \bbN$, be a family of
admissible ordered sets. Let $\calL$ be the Krull--Schmidt
category, whose indecomposable objects are
\begin{itemize}

\item
$X_\gamma$, $\gamma \in I_p'$, $p \in [0, r + 1]$,

\item
$X_{\max I_p}'$, $X_{\max I_p}''$, $p \in [0, r]$,

\item
$Y_\gamma$, $\gamma \in I_0'$,

\item
$Z$.

\end{itemize}
If $U$ is an indecomposable object of $\calL$, then
\[
\dim_k |U| =
\begin{cases}
2 & U = X_{\min I_1}, \\ %
1 & \text{otherwise}.
\end{cases}
\]
If $U$ and $V$ are indecomposable objects of $\calL$, then $\dim_k
\Hom_{\calL} (U, V) \leq 2$, $\Hom_{\calL} (U, V) \neq 0$ if and
only if one of the following conditions holds:
\begin{itemize}

\item
$U = X_{\gamma'}$, $V = X_{\gamma''}$, $\gamma' \in I_p'$,
$\gamma'' \in I_q'$, $(p, \gamma') \leq (q, \gamma'')$,

\item
$U = X_\gamma$, $V = X_{\max I_q}'$, $\gamma \in I_p'$, $p \leq
q$,

\item
$U = X_\gamma$, $V = X_{\max I_q}''$, $\gamma \in I_p'$, $p \leq
q$,

\item
$U = X_{\gamma'}$, $V = Y_{\gamma''}$, $\gamma', \gamma'' \in
I_0'$, $\gamma' \leq \gamma''$,

\item
$U = X_\gamma$, $V = Z$, $\gamma \in I_0'$,

\item
$U = X_{\max I_p}'$, $V = X_\gamma$, $\gamma \in I_q'$, $p < q$,

\item
$U = X_{\max I_p}'$, $V = X_{\max I_q}'$, $p \leq q$,

\item
$U = X_{\max I_p}'$, $V = X_{\max I_q}''$, $p < q$,

\item
$U = X_{\max I_p}''$, $V = X_\gamma$, $\gamma \in I_q'$, $p < q$,

\item
$U = X_{\max I_p}''$, $V = X_{\max I_q}'$, $p < q$,

\item
$U = X_{\max I_p}''$, $V = X_{\max I_q}''$, $p \leq q$,

\item
$U = X_{\max I_0}''$, $V = Z$,

\item
$U = Y_\gamma$, $V = X_{\min I_1}$,

\item
$U = Y_{\gamma'}$, $V = Y_{\gamma''}$, $\gamma' \leq \gamma''$,

\item
$U = Y_\gamma$, $V = Z$,

\item
$U = Z$, $V = Z$,

\end{itemize}
and $\dim_k \Hom_{\calL} (U, V) = 2$ if and only if $U = X_\gamma$,
$\gamma \in I_0'$, $V = X_{\min I_1}$. By $\bbL_{I_0, \ldots, I_{r +
1}}$ we denote the vector space category $(\calL, {|-|})$, where
${|-|} : \calL \to \mod k$ is the forgetful functor. We refer the
reader to~\cite{BobSk1}*{Section~1} for pictures presenting vector
space categories of the above type, and in particular explaining how
the forgetful functor ${|-|}$ is defined on $\Hom_{\calL} (X_\gamma,
X_{\min I_1})$ for $\gamma \in I_0'$.

\subsection{} \label{propLIr}
We describe the indecomposable objects and the Auslander--Rei\-ten
sequences in $\calU (\bbL_{I_0, \ldots, I_{r + 1}})$. We refer
to~\cite{BobSk1} for definitions of the objects listed in the below
proposition and its proof.

\begin{prop*}
Let $I_0$, \ldots, $I_{r + 1}$, $r \in \bbN$, be admissible
ordered sets. Put
\[
I_p'' =
\begin{cases}
I_1' \setminus \{ \min I_1 \} & p = 1, \\ %
I_p' & p \in [2, r + 1],
\end{cases}
\]
\begin{enumerate}

\item
Objects
\begin{align*}
& M_{(-1, \max I_0), (0, \gamma)} = X_\gamma, \, \gamma \in I_0',
\\ %
& M_{(0, \gamma'), (0, \gamma'')} = \ol{Y_{\gamma'} X_{\gamma''}},
\, \gamma', \gamma'' \in I_0', \, \gamma' < \gamma'',
\\ %
& M_{(n - 1, \max I_0), (n, \gamma)} = \ol{Y_\gamma X_{\max I_0}'
X_{\max I_0}'' X_{\min I_1}^{2 n - 1}}^{2 n}, \, \gamma \in I_0',
\, n \in \bbN,
\\ %
& M_{(n, \gamma), (n, \max I_0)} = \ol{Y_\gamma X_{\max I_0}'
X_{\max I_0}'' X_{\min I_1}^{2 n}}^{2 n + 1}, \, \gamma \in I_0',
\, n \in \bbN_0,
\\ %
&  M_{(n, \gamma''), (n + 1, \gamma')} = \ol{Y_{\gamma''}
Y_{\gamma'} X_{\max I_0}' X_{\max I_0}'' X_{\min I_1}^{2 n}}^{2 n
+ 2},
\\ %
& \qquad \gamma', \gamma'' \in I_0', \gamma' < \gamma'', \, n \in
\bbN_0,
\\ %
& M_{(n, \gamma'), (n, \gamma'')} = \ol{Y_{\gamma'} Y_{\gamma''}
X_{\max I_0}' X_{\max I_0}'' X_{\min I_1}^{2 n - 1}}^{2 n + 1},
\\ %
& \qquad  \gamma', \gamma'' \in I_0', \gamma' < \gamma'', \, n \in
\bbN,
\\ %
& M_{(n, \gamma), (n, \gamma)}' = \ol{Y_\gamma X_{\min I_1}^n}^n,
\gamma \in I_0', \, n \in \bbN_0,
\\ %
& M_{(n, \max I_0), (n, \max I_0)}' = \ol{X_{\min I_1}^{n + 1}}^n,
\, n \in \bbN_0,
\\ %
& M_{(0, \gamma), (0, \gamma)}'' = \ol{X_\gamma}, \, \gamma \in
I_0',
\\ %
& M_{(n, \max I_0), (n, \max I_0)}'' = \ol{X_{\max I_0}' X_{\max
I_0}'' X_{\min I_1}^n}^{n + 1}, \, n \in \bbN_0,
\\ %
& M_{(n, \gamma), (n, \gamma)}'' = \ol{Y_\gamma X_{\max I_0}'
X_{\max I_0}'' X_{\min I_1}^{n - 1}}^{n + 1}, \, \gamma \in I_0',
\, n \in \bbN,
\\ %
& M_{(n - 1, \max I_0), (n, \max I_0)}' = \ol{X_{\max I_0}'
X_{\min I_1}^n}^n, \, n \in \bbN_0,
\\ %
& M_{(n, \gamma), (n + 1, \gamma)}' = \ol{Y_\gamma X_{\max I_0}'
X_{\min I_1}^n}^{n + 1}, \, \gamma \in I_0', \, n \in \bbN_0,
\\ %
& M_{(n - 1, \max I_0), (n, \max I_0)}'' = \ol{X_{\max I_0}''
X_{\min I_1}^n}^n, \, n \in \bbN_0,
\\ %
& M_{(n, \gamma), (n + 1, \gamma)}'' = \ol{Y_\gamma X_{\max I_0}''
X_{\min I_1}^n}^{n + 1}, \, \gamma \in I_0', \, n \in \bbN_0,
\\ %
& R_n^\lambda = \ol{X_{\min I_1}^n}^n (\lambda), \, \lambda \in
k^*, \, \lambda \neq 1, \, n \in \bbN,
\\ %
& R_{2 n - 1}^1 = \ol{X_{\max I_0}'' X_{\min I_1}^{n - 1}}^n, \, n
\in \bbN,
\\ %
& R_{2 n}^1 = \ol{X_{\min I_1}^n}^n (1), \, n \geq 1,
\\ %
& R_{2 n - 1, 0}^\infty = \ol{X_{\max I_0}' X_{\min I_1}^{n -
1}}^n, \, n \in \bbN,
\\ %
& R_{2 n - 1, 1}^\infty = \ol{X_{\min I_1}^{n - 1} Z}^{n - 1}, \,
n \in \bbN,
\\ %
& R_{2 n, 0}^\infty = \ol{X_{\min I_1}^n}^n (\infty), \, n \in
\bbN,
\\ %
& R_{2 n, 1}^\infty = \ol{X_{\max I_0}' X_{\min I_1}^{n - 1} Z}^n,
\, n \in \bbN,
\\ %
& S_{p, (n - 1, \max I_0), (m - 1, \max I_0)} = \ol{X_{\min I_1}^n
X_{\max I_p}' X_{\max I_p}'' X_{\min I_1}^m}^{n + m}, \\ %
& \qquad p \in [1, r], \, n, m \in \bbN_0, \, n < m,
\\ %
& S_{p, (n, \gamma), (m - 1, \max I_0)} = \ol{Y_\gamma X_{\min
I_1}^n X_{\max I_p}' X_{\max I_p}'' X_{\min I_1}^m}^{n + m + 1},
\\ %
& \qquad p \in [1, r], \, \gamma \in I_0', \, n, m \in \bbN_0, \,
n < m,
\\ %
& S_{p, (n - 1, \max I_0), (m, \gamma)} = \ol{X_{\min I_1}^n
X_{\max I_p}' X_{\max I_p}'' X_{\min I_1}^m Y_\gamma}^{n + m + 1},
\\ %
& \qquad p \in [1, r], \, \gamma \in I_0', \, n, m \in \bbN_0, \,
n \leq m,
\\ %
& S_{p, (n, \gamma'), (m, \gamma'')} = \ol{Y_{\gamma''} X_{\min
I_1}^n X_{\max I_p}' X_{\max I_p}'' X_{\min I_1}^m
Y_{\gamma'}}^{n + m + 2}, \\ %
& \qquad p \in [1, r], \, \gamma, \gamma' \in I_0', \, n, m \in
\bbN_0, \, (n, \gamma') < (m, \gamma''),
\\ %
& S_{p, (n - 1, \max I_0), (n - 1, \max I_0)}' = \ol{X_{\min
I_1}^n X_{\max I_p}'}^n, \, p \in [1, r], \, n \in \bbN_0,
\\ %
& S_{p, (n, \gamma), (n, \gamma)}' = \ol{Y_\gamma X_{\min I_1}^n
X_{\max I_p}'}^{n + 1}, \, p \in [1, r], \, \gamma \in I_0', \, n
\in \bbN_0,
\\ %
& S_{p, (n - 1, \max I_0), (n - 1, \max I_0)}'' = \ol{X_{\min
I_1}^n X_{\max I_p}''}^n, \, p \in [1, r], \, n \in \bbN_0,
\\ %
& S_{p, (n, \gamma), (n, \gamma)}'' = \ol{Y_\gamma X_{\min I_1}^n
X_{\max I_p}''}^{n + 1}, \, p \in [1, r], \, \gamma \in I_0', \, n
\in \bbN_0,
\\ %
& T_{p, \gamma, (m - 1, \max I_0)} = \ol{X_{\min I_1}^m
X_\gamma}^m, \, p \in [1, r + 1], \, \gamma \in I_p'', \, m \in
\bbN_0,
\\ %
& T_{p, \gamma', (m, \gamma'')} = \ol{Y_{\gamma''} X_{\min I_1}^m
X_{\gamma'}}^{m + 1}, \, p \in [1, r + 1], \, \gamma' \in I_p'',
\, \gamma'' \in I_0', \, m \in \bbN_0,
\\ %
& T_{r + 1, \max I_{r + 1}, (m - 1, \max I_0)} = \ol{X_{\min
I_1}^m}^m (0), \, m \in \bbN,
\\ %
& T_{r + 1, \max I_{r + 1}, (m, \gamma)} = \ol{Y_\gamma X_{\min
I_1}^m}^{m + 1}, \, \gamma \in I_0', \, m \in \bbN_0,
\\ %
& U_{p, 2 n, (m - 1, \max I_0)} = \ol{X_{\min I_1}^m X_{\max I_p}'
X_{\max I_p}'' X_{\min I_1}^n}^{n + m + 1},
\\ %
& \qquad p \in [1, r], \, n, m \in \bbN_0,
\\ %
& U_{p, 2 n, (m, \gamma)} = \ol{Y_\gamma X_{\min I_1}^m X_{\max
I_p}' X_{\max I_p}'' X_{\min I_1}^n}^{n + m + 2},
\\ %
& \qquad p \in [1, r], \, \gamma \in I_0', \, n, m \in \bbN_0,
\\ %
& U_{p, 2 n + 1, (m - 1, \max I_0)} = \ol{X_{\min I_1}^m X_{\max
I_p}' X_{\max I_p}'' X_{\min I_1}^n Z}^{n + m + 1},
\\ %
& \qquad p \in [1, r], \, n, m \in \bbN_0,
\\ %
& U_{p, 2 n + 1, (m, \gamma)} = \ol{Y_\gamma X_{\min I_1}^m
X_{\max I_p}' X_{\max I_p}'' X_{\min I_1}^n Z}^{n + m + 2},
\\ %
& \qquad p \in [1, r], \, \gamma \in I_0', \, n, m \in \bbN_0,
\\ %
& V_{p, 2 n, \gamma} = \ol{X_{\min I_1}^n X_\gamma}^{n + 1}, \, p
\in [1, r + 1], \, \gamma \in I_p'', \, n \in \bbN_0,
\\ %
& V_{p, 2 n + 1, \gamma} = \ol{X_{\min I_1}^n X_\gamma Z}^{n + 1},
\, p \in [1, r + 1], \, \gamma \in I_p'', \, n \in \bbN_0,
\\ %
& V_{r + 1, 2 n, \max I_{r + 1}} = \ol{X_{\min I_1}^n}^{n + 1}, \,
n \in \bbN_0,
\\ %
& V_{r + 1, 2 n + 1, \max I_{r + 1}} = \ol{X_{\min I_1}^n Z}^{n +
1}, \, n \in \bbN_0,
\\ %
& W_{p, 2 n, 2 m} = \ol{X_{\min I_1}^n X_{\max I_p}' X_{\max
I_p}'' X_{\min I_1}^m}^{n + m + 2},
\\ %
& \qquad p \in [1, r], \, n, m \in \bbN_0, \, m < n,
\\ %
& W_{p, 2 n + 1, 2 m} = \ol{Z X_{\min I_1}^n X_{\max I_p}' X_{\max
I_p}'' X_{\min I_1}^m}^{n + m + 2},
\\ %
& \qquad p \in [1, r], \, n, m \in \bbN_0, \, m \leq n,
\\ %
& W_{p, 2 n, 2 m + 1} = \ol{X_{\min I_1}^n X_{\max I_p}' X_{\max
I_p}'' X_{\min I_1}^m Z}^{n + m + 2},
\\ %
& \qquad p \in [1, r], \, n, m \in \bbN_0, \, m < n,
\\ %
& W_{p, 2 n + 1, 2 m + 1} = \ol{Z X_{\min I_1}^n X_{\max I_p}'
X_{\max I_p}'' X_{\min I_1}^m Z}^{n + m + 2},
\\ %
& \qquad p \in [1, r], \, n, m \in \bbN_0, \, m < n,
\\ %
& W_{p, 2 n, 2 n}' = \ol{X_{\min I_1}^n X_{\max I_p}'}^{n + 1}, \,
p \in [1, r], \, n \in \bbN_0,
\\ %
& W_{p, 2 n + 1, 2 n + 1}' = \ol{X_{\min I_1}^n X_{\max I_p}'
Z}^{n + 1}, \, p \in [1, r], \, n \in \bbN_0,
\\ %
& W_{p, 2 n, 2 n}'' = \ol{X_{\min I_1}^n X_{\max I_p}''}^{n + 1},
\, p \in [1, r], \, n \in \bbN_0,
\\ %
& W_{p, 2 n + 1, 2 n + 1}'' = \ol{X_{\min I_1}^n X_{\max I_p}''
Z}^{n + 1}, \, p \in [1, r], \, n \in \bbN_0,
\end{align*}
form a complete list of indecomposable objects in $\calU
(\bbL_{I_0, \ldots, I_{r + 1}})$.

\item
Sequences
\begin{align*}
& 0 \to M_{\gamma', \gamma''} \to M_{\gamma'_+, \gamma''} \oplus
M_{\gamma', \gamma''_+} \to M_{\gamma'_+, \gamma''_+} \to 0,
\\ %
& \qquad \gamma', \gamma'' \in \bbZ \times I_0, \, (-1, \max I_0)
\leq \gamma' < \gamma'' < (\gamma')^+,
\\ %
& 0 \to M_{\gamma, \gamma}' \to M_{\gamma, \gamma_+} \to
M_{\gamma_+, \gamma_+}'' \to 0, \, \gamma \in \bbN_0 \times I_0,
\\ %
& 0 \to M_{\gamma, \gamma}'' \to M_{\gamma, \gamma_+} \to
M_{\gamma_+, \gamma_+}' \to 0, \, \gamma \in \bbN_0 \times I_0,
\\ %
& 0 \to M_{\gamma, \gamma^+}' \to M_{\gamma_+, \gamma^+} \to
M_{\gamma_+, \gamma_+^+}'' \to 0, \, \gamma \in \bbZ \times I_0,
\, (-1, \max I_0) \leq \gamma,
\\ %
& 0 \to M_{\gamma, \gamma^+}'' \to M_{\gamma_+, \gamma^+} \to
M_{\gamma_+, \gamma_+^+}' \to 0, \, \gamma \in \bbZ \times I_0, \,
(-1, \max I_0) \leq \gamma,
\\ %
& 0 \to R_n^\lambda \to R_{n + 1}^\lambda \oplus R_{n - 1}^\lambda
\to R_n^\lambda \to 0, \, \lambda \in k^*, \, \lambda \neq 1, \, n
\in \bbN,
\\ %
& 0 \to R_{n + 1}^1 \to R_{n + 2}^1 \oplus R_{n - 1}^1 \to R_n^1
\to 0, \, n \in \bbN,
\\ %
& 0 \to R_{n, i}^\infty \to R_{n + 1, i}^\infty \oplus R_{n - 1, i
+ 1}^\infty \to R_{n, i + 1}^\infty \to 0, \, i \in \bbZ_2, \, n
\in \bbN,
\\ %
& 0 \to S_{p, \gamma', \gamma''} \to S_{p, \gamma'_+, \gamma''}
\oplus S_{p, \gamma', \gamma''_+} \to S_{p, \gamma'_+, \gamma''_+}
\to 0,
\\ %
& \qquad p \in [1, r], \, \gamma', \gamma'' \in \bbZ \times I_0,
\, (-1, \max I_0) \leq \gamma' < \gamma'',
\\ %
& 0 \to S_{p, \gamma, \gamma}' \to S_{p, \gamma, \gamma_+} \to
S_{p, \gamma_+, \gamma_+}'' \to 0,
\\ %
& \qquad p \in [1, r], \, \gamma \in \bbZ \times I_0, \, (-1, \max
I_0) \leq \gamma,
\\ %
& 0 \to S_{p, \gamma, \gamma}'' \to S_{p, \gamma, \gamma_+} \to
S_{p, \gamma_+, \gamma_+}' \to 0,
\\ %
& \qquad p \in [1, r], \, \gamma \in \bbZ \times I_0, \, (-1, \max
I_0) \leq \gamma,
\\ %
& 0 \to T_{r + 1, \max I_{r + 1}, (0, \max I_0')} \to T_{r + 1,
\max I_{r + 1}, (0 \max I_0)}
\\ %
& \qquad \to T_{1, \min I_1'', (0, \max I_0)} \to 0,
\\ %
& 0 \to T_{p, \gamma', \gamma''} \to T_{p, \gamma'_+, \gamma''}
\oplus T_{p, \gamma', \gamma''_+} \to T_{p, \gamma'_+, \gamma''_+}
\to 0, \, p \in [1, r + 1],
\\ %
& \qquad \gamma' \in (I_p'')_-,\, \gamma'' \in \bbZ \times I_0, \,
(-1, \max I_0) \leq \gamma'',
\\ %
& 0 \to U_{p, n, \gamma} \to U_{p, n, \gamma_+} \oplus U_{p, n -
1, \gamma} \to U_{p, n - 1, \gamma_+} \to 0,
\\ %
& \qquad p \in [1, r], \, n \in \bbN, \, \gamma \in \bbN_0 \times
I_0,
\\ %
& 0 \to V_{p, n, \gamma} \to V_{p, n, \gamma_+} \oplus V_{p, n -
1, \gamma} \to V_{p, n - 1, \gamma_+} \to 0,
\\ %
& \qquad p \in [1, r + 1], \, n \in \bbN, \, \gamma \in (I_p'')_-,
\\ %
& 0 \to W_{p, n, m} \to W_{p, n - 1, m} \oplus W_{p, n, m - 1} \to
W_{p, n - 1, m - 1} \to 0,
\\ %
& \qquad p \in [1, r], \, n, m \in \bbN, m < n,
\\ %
& 0 \to W_{p, n, n}' \to W_{p, n, n - 1} \to W_{p, n - 1, n - 1}''
\to 0, \, p \in [1, r], \, n \in \bbN,
\\ %
& 0 \to W_{p, n, n}'' \to W_{p, n, n - 1} \to W_{p, n - 1, n - 1}'
\to 0, \, p \in [1, r], \, n \in \bbN,
\end{align*}
form a complete list of Auslander--Reiten sequences in the
category $\calU (\bbL_{I_0, \ldots, I_{r + 1}})$, where
\begin{align*}
& M_{\gamma, \gamma} = M_{\gamma, \gamma}' \oplus M_{\gamma,
\gamma}'', \, \gamma \in \bbN_0 \times I_0,
\\ %
& M_{\gamma, \gamma^+} = M_{\gamma, \gamma^+}' \oplus M_{\gamma,
\gamma^+}'', \, \gamma \in \bbZ \times I_0, \, (-1, \max I_0) \leq
\gamma,
\\ %
& R_0^\lambda = 0, \, \lambda \in k^*,
\\ %
& R_{0, i}^\infty = 0, \, i \in \bbZ_2,
\\ %
& S_{\gamma, \gamma} = S_{\gamma, \gamma}' \oplus S_{\gamma,
\gamma}'', \, \gamma \in \bbN_0 \times I_0,
\\ 
& T_{1, \min (I_1'')_-, \gamma} = T_{r + 1, \max I_{r + 1},
\gamma^+}, \, \gamma \in \bbZ \times I_0, \, (-1, \max I_0) \leq
\gamma,
\\ %
& T_{p, \min (I_p'')_-, \gamma} = U_{p - 1, 0, \gamma}, \, p \in
[2, r + 1], \, \gamma \in \bbZ \times I_0, \, (-1, \max I_0) \leq
\gamma,
\\ %
& T_{p, \max I_p, \gamma} = S_{p, (-1, \max I_0), \gamma}, \, p
\in [1, r], \, \gamma \in \bbZ \times I_0, \, (-1, \max I_0) \leq
\gamma,
\\ %
& T_{r + 1, \max I_{r + 1}, (-1, \max I_0)} = 0,
\\ %
& V_{1, n, \min (I_1'')_-} = V_{r + 1, n + 2, \max I_{r + 1}}, \,
n \in \bbN_0,
\\ %
& V_{p, n, \min (I_p'')_-} = W_{p - 1, n, 0}, \, p \in [2, r + 1],
\, n \in \bbN_0,
\\ %
& V_{p, n, \max I_p} = U_{p, n, (-1, \max I_0)}, \, p \in [1, r],
\, n \in \bbN_0,
\\ %
& W_{p, n, n} = W_{p, n, n}' \oplus W_{p, n, n}'', \, p \in [1,
r], \, n \in \bbN_0.
\end{align*}

\end{enumerate}
\end{prop*}

\section{Proof of the main result} \label{sectproof}

In this section we present the proof of the main theorem of the
paper.

\subsection{}
Let $A$ be an algebra and let $R$ be an $A$-module. By $A [R]$ we
denote the one-point extension of $A$ by $R$ defined as
\[
\begin{bmatrix}
A & R \\ 0 & k
\end{bmatrix}.
\]
The category of $A [R]$-modules is equivalent to the category of
triples $(V_0, V_1, \gamma_V)$, with $V_0 \in \mod A$, $V_1 \in
\mod k$ and $\gamma_V : V_1 \to \Hom_A (R, V_0)$ is a $k$-linear
map (see \cite{Ri2}*{2.5(8)}).

Let $\Hom (R, \mod A)$ be the vector space category $(\calK,
{|-|})$, where $\calK = \mod A / \Ker \Hom_A (R, -)$ and ${|-|} :
\calK \to \mod k$ is the functor induced by $\Hom_A (R, -)$. It
follows from the above remark that we may view the objects of
$\calU (\Hom (R, \mod A))$ as objects of $\mod A [R]$.
Consequently, if $X$ is an indecomposable $A [R]$-module then
either $X \in \mod A$ or $X \in \calU (\Hom (R, \mod A))$.
Moreover, each Auslander--Reiten sequence in $\mod A [R]$ is
either of the form
\begin{multline*}
0 \to (X, \Hom_A (R, X), \Hom_A (R, \Id_X))
\\ %
\to (Y, \Hom_A (R, X), \Hom_A (R, f)) \to (Z, 0, 0) \to 0
\end{multline*}
for an Auslander--Reiten sequence $0 \to X \xrightarrow{f} Y \to Z
\to 0$ in $\mod A$, or comes from an Auslander--Reiten sequence in
$\calU (\Hom (R, \mod A))$.

\subsection{}
From now on we assume that $(p, q, S, T)$ is a fixed defining
system. We also use notation introduced in Section~\ref{mainres}.

A vertex $x$ of $Q$ is called admissible if one of the following
possibilities holds:
\begin{itemize}

\item
$x = x_{i, j}$, $i \in [1, |p|]$, $j \in [2, p_i + |T_i|]$, $j -
1, j, j + 1 \not \in S_i$,

\item
$x = z_{i, j}$, $i \in [1, |p|]$, $j \in S_i \cap [T_{i, |T_i|} +
2, p_i + |T_i|]$.

\end{itemize}
For an admissible vertex $x$ of $Q$ we define a defining system
$(p, q, S^x, T^x)$ by:
\begin{itemize}

\item
if $x = x_{i_0, j_0}$, then
\[
S_i^x =
\begin{cases}
S_{i_0} \cup \{ j_0 \} & i = i_0,
\\ %
S_i & i \neq i_0,
\end{cases}
\text{ and } T^x = T,
\]

\item
if $x = z_{i_0, j_0}$, then
\[
S^x = S \text{ and } T_i^x =
\begin{cases}
T_{i_0} \cup \{ j_0 \} & i = i_0,
\\ %
T_i & i \neq i_0.
\end{cases}
\]

\end{itemize}
A defining system $(p, q, S, T)$ if called fundamental if $S_i =
\varnothing = T_i$ for all $i$. The following observation allows
us to perform inductive proofs: each defining system is an
iterated extension of a fundamental one by admissible vertices.

\subsection{} \label{sectlemm}
For $x \in Q_0$, $x = x_{i, j}$, let $X_x = M (\mu_x)$, $I_x = M
(\omega_x)$ and $R_x = M (\alpha_x \mu_x)$. Similarly, if $x =
z_{i, j}$, then $X_x = M (\gamma_y \mu_y)$ and $R_x = N (\mu_y,
\omega_y)$, where $y = x_{i, j}$.

For a vertex $x$ of $Q$ let $\calC_x$ denote the set of all strings
terminating at $x$ ordered by the relation introduced
in~\ref{sectord}. Recall that $\calC_x' = \calC_x \setminus \{
\omega_x \}$. We prove the main theorem inductively together with
the following series of lemmas.

\begin{lemm} \label{lemmR}
Let $x$ be an admissible vertex of $Q$.
\begin{enumerate}

\item
If $x = x_{i_0, j_0}$, then the assignment
\begin{align*}
X_C & \mapsto M (\alpha_x C), \, C \in \calC_x,
\\ %
Y_C & \mapsto M (C), \, C \in \calC_x,
\end{align*}
induces an equivalence between $\bbL_{\calC_x}$ and $\Hom (R_x,
\mod A)$.

\item
If $x = z_{i_0, j_0}$, let $\{ j_1 < \cdots < j_r \} = S_{i_0}
\cap [j_0 + 1, p_{i_0} + |T_{i_0}|]$ and $j_{r + 1} = p_{i_0} +
|T_{i_0}| + 1$. The assignment
\begin{align*}
X_C & \mapsto N (C, \omega_{x_{i_0, j_p}}), \, C \in
\calC_{x_{i_0, j_p}}', \, p \in [0, r],
\\ %
X_{j_p} & \mapsto M (\gamma_{i_0, j_p} \omega_{x_{i_0, j_p}}), \,
p \in [0, r],
\\ %
X_j & \mapsto M (\omega_{x_{i_0, j}}), \, j \in [j_p + 1, \ldots,
j_{p + 1} - 1], \, p \in [0, r],
\\ %
X_{\omega_{x_{i_0, j_p}}}' & \mapsto M (\omega_{x_{i_0, j_p}}), \,
p \in [0, r],
\\ %
X_{\omega_{x_{i_0, j_p}}}'' & \mapsto N (\omega_{x_{i_0, j_p}}),
\, p \in [0, r],
\\ %
Y_C & \mapsto M (\gamma_{i_0, j_0} C), \, C \in \calC_{x_{i_0,
j_0}}',
\\ %
Z & \mapsto M (x),
\end{align*}
induces an equivalence between
\[
\bbL_{\calC_{x_{i_0, j_0}}, [j_0, j_1 - 1] + \calC_{x_{i_0, j_1}},
\ldots, [j_{r - 1}, j_r - 1] + \calC_{x_{i_0, j_r}}, [j_r, j_{r +
1}]} \text{ and } \Hom (R_x, \mod A).
\]

\end{enumerate}
\end{lemm}

\begin{lemm} \label{lemmX}
Let $x$ be an admissible vertex of $Q$. The assignment
\[
X_C \mapsto M (C), \, C \in \calC_x,
\]
induces an equivalence between $\bbK_{\calC_x}$ and $\Hom (X_x,
\mod A)$.
\end{lemm}

\begin{lemm} \label{lemmI}
Let $x = x_{i_0, j_0}$ be such that $j_0 \in [T_{i_0, |T_{i_0}|} +
1, p_{i_0} + |T_{i_0}|] \setminus S_{i_0}$. Let $\{ j_1 < \cdots <
j_r \} = S_{i_0} \cap [j_0 + 1, p_{i_0} + |T_{i_0}|]$ and $j_{r +
1} = p_{i_0} + |T_{i_0}| + 1$. The assignment
\begin{align*}
X_C & \mapsto N (C, \omega_{x_{i_0, j_p}}), \, C \in
\calC_{x_{i_0, j_p}}', \, p \in [1, r],
\\ %
X_{j_0} & \mapsto M (\omega_{x_{i_0, j_0}}),
\\ %
X_{j_p} & \mapsto M (\gamma_{i_0, j_p} \omega_{x_{i_0, j_p}}), \,
p \in [1, r],
\\ %
X_j & \mapsto M (\omega_{x_{i_0, j}}), \, j \in [j_p + 1, j_{p
+ 1} - 1], \, p \in [0, r],
\\ %
X_{\omega_{x_{i_0, j_p}}}' & \mapsto M (\omega_{x_{i_0, j_p}}), \,
p \in [1, r],
\\ %
X_{\omega_{x_{i_0, j_p}}}'' & \mapsto N (\omega_{x_{i_0, j_p}}),
\, p \in [1, r],
\end{align*}
induces an equivalence between
\[
\bbK_{[j_0, j_1 - 1] + \calC_{x_{i_0, j_1}}, \ldots, [j_{r - 1},
j_r - 1] + \calC_{x_{i_0, j_r}}, [j_r, j_{r + 1} - 1]} \text{ and
} \Hom (I_x, \mod A).
\]
\end{lemm}

\subsection{}
If $(p, q, S, T)$ is a fundamental defining system, then
Theorem~\ref{maintheo} and Lemmas~\ref{sectlemm} are easy
exercises in the representation theory of a hereditary algebra of
type $\tilde{\bbA}_{p, q}$.

From now on we assume that we Theorem~\ref{maintheo} and
Lemmas~\ref{sectlemm} have been proved for $(p, q, S, T)$. Let $x$
be an admissible vertex of $Q$. We will show that
Theorem~\ref{maintheo} and Lemmas~\ref{sectlemm} hold for $(p, q,
S^x, T^x)$.

By $Q^x$ (respectively, $A^x$) we will denote the quiver (algebra)
associated with $(p, q, S^x, T^x)$. We also define $R_{x'}^x$,
$X_{x'}^x$ and $I_{x'}^x$ in the analogous way as the
corresponding modules for $(p, q, S, T)$.

\subsection{}
Assume first that $x = x_{i_0, j_0}$. Let $\gamma = \gamma_{i_0,
j_0}$ be the new arrow of $Q^x$ and $z = z_{i_0, j_0}$ be the new
vertex of $Q^x$. Theorem~\ref{maintheo} for $(p, q, S^x, T^x)$
follows from the induction hypothesis (Theorem~\ref{maintheo} and
Lemma~\ref{lemmR} for $(p, q, S, T)$), Proposition~\ref{subspaceone}
and the following isomorphisms
\begin{align*}
& \ol{M (\alpha_x C)} \simeq N (C), \, C \in \calC_x,
\\ %
& \ol{M (C)} \simeq M (\gamma C), \, C \in \calC_x,
\\ %
& \ol{M (C') M (\alpha_x C'')} \simeq N (C', C''), \, C', C'' \in
\calC_x, \, C' < C'',
\\ %
& \ol{0} \simeq M (z).
\end{align*}

\subsection{}
Now we prove Lemma~\ref{lemmR} for $(p, q, S^x, T^x)$. Let $x'$ be
an admissible vertex of $Q^x$. Then either $x'$ is an admissible
vertex of $Q$ or $x' = z$. In the first case there are still two
possibilities: either $x' = x_{i, j}$ or $x' = z_{i, j}$.

Consider first the case $x' = x_{i, j}$. Then either $i \neq i_0$ or
$i = i_0$ and $|j - j_0| > 1$, hence it is easily seen that in this
case we also have $R_{x'}^x = R_{x'}$ and $\Hom (R_{x'}^x, \mod A^x)
= \Hom (R_{x'}, \mod A)$, thus the claim follows.

Let now $x' = z_{i, j}$ for $(i, j) \neq (i_0, j_0)$. In this case
also $R_{x'}^x = R_{x'}$. Moreover, if $i \neq i_0$ or $i = i_0$ and
$j_0 < j$, then $\Hom (R_{x'}^x, \mod A^x) = \Hom (R_{x'}, \mod A)$.
If $j < j_0$, then the claim about $\Hom (R_{x'}^x, \mod A^x)$
follows by observing that its indecomposable objects are the
indecomposable objects of $\Hom (R_{x'}, \mod A)$ and
\[
\ol{M (C) M (\alpha_x \omega_x)}, \, C \in \calC_x', \, \ol{M
(\alpha_x \omega_x)}, \, \ol{M (\omega_x)}.
\]

Finally, let $x' = z$. Then $R_{x'}^x = \ol{I_x X_x}$ and the
claim follows from Lemmas~\ref{lemmX} and \ref{lemmI}.

\subsection{}
In order to show Lemma~\ref{lemmX} we have to consider the cases
analogous to the ones considered above. If $x' = x_{i, j}$ or $x' =
z_{i, j}$, $x' \neq z$, then $X_{x'}^x = X_{x'}$ and $\Hom
(X_{x'}^x, \mod A^x) = \Hom (X_{x'}, \mod A)$. Thus it remains to
consider the case $x' = z$. In this case $X_{x'}^x = \ol{X_x}$ and
the description of $\Hom (X_{x'}^x, \mod A^x)$ follows easily from
the description of $\Hom (X_x, \mod A)$.

\subsection{}
It remains to show Lemma~\ref{lemmI}. Let $x' = x_{i, j}$ be the
vertex of $Q^x$ satisfying the hypothesis of Lemma~\ref{lemmI}. We
have $I_{x'}^x = I_{x'}$. If $i \neq i_0$ or $i = i_0$ and $j_0 <
j$, then also $\Hom (I_{x'}^x, \mod A^x) = \Hom (I_{x'}, \mod A)$.
If $j < j_0$, then we have observed that indecomposable objects of
$\Hom (I_{x'}^x, \mod A^x)$ are the indecomposable objects of $\Hom
(I_{x'}, \mod A)$ and
\[
\ol{M (C) M (\alpha_x \omega_x)}, \, C \in \calC_x', \, \ol{M
(\alpha_x \omega_x)}, \, \ol{M (\omega_x)}.
\]

\subsection{}
Assume now that $x = z_{i_0, j_0}$. Let $\{ j_1 < \cdots < j_r \}
= S_{i_0} \cap [j + 1, p_{i_0} + |T_{i_0}|]$ and $j_{r + 1} =
p_{i_0} + |T_{i_0}| + 1$. Put
\begin{align*}
\calC & = \calC_{x_{i_0, j_0}}, &
\calC_p & = \calC_{x_{i_0, j_p}}, \, p \in [1, r],
\\ %
\gamma & = \gamma_{i_0, j_0}, &
\gamma_p & = \gamma_{i_0, j_p}, \, p \in [1, r],
\\ %
\omega & = \omega_{x_{i_0, j_0}}, &
\omega_j & = \omega_{x_{i_0, j}}, \, j \in [j + 1, p_{i_0} +
|T_{i_0}|],
\\ %
\alpha & = \alpha_{i_0, p_{i_0} + |T_{i_0}| + 1}, &
\xi & = \xi_{i_0, |T_{i_0}| + 1},
\\ %
\intertext{and} %
B & = \omega \alpha \xi \gamma, &
B_j & = \omega_j \alpha \xi \gamma, \, j \in [j + 1, p_{i_0} +
|T_{i_0}|].
\end{align*}
Finally, let $z = x_{i_0, j_{r + 1}}$.

\subsection{}
In this case Theorem~\ref{maintheo} for $(p, q, S^x, T^x)$ follows
from the induction hypothesis, Proposition~\ref{propLIr} and the
following isomorphisms
\begin{align*}
& \ol{M (\gamma C') N (C'', \omega)} \simeq N (C'', B C'), \, C',
C'' \in \calC', \, C' < C'',
\\ %
& \ol{M (\gamma C) M (\omega) N (\omega) M (\gamma \omega)^{2 n -
1}}^{2 n} \simeq N (B^n C, B^n \omega), \, C \in \calC', \, n \in
\bbN,
\\ %
& \ol{M (\gamma C) M (\omega) N (\omega) M (\gamma \omega)^{2 n}}^{2
n + 1} \simeq N (B^n \omega, B^{n + 1} C), \, C \in \calC', \, n \in
\bbN,
\\ %
& \ol{M (\gamma C'') M (\gamma C') M (\omega) N (\omega) M (\gamma
\omega)^{2 n}}^{2 n + 2} \simeq N (B^{n + 1} C', B^{n + 1} C''),
\\ %
& \qquad C', C'' \in \calC', \, C' < C'', \, n \in \bbN_0,
\\ %
& \ol{M (\gamma C') M (\gamma C'') M (\omega) N (\omega) M (\gamma
\omega)^{2 n - 1}}^{2 n + 1} \simeq N (B^n C'', B^{n + 1} C'),
\\ %
& \qquad C', C'' \in \calC', \, C' < C'', \, n \in \bbN,
\\ %
& \ol{M (\gamma C) M (\gamma \omega)^n}^n \simeq M (\gamma B^n C),
\, C \in \calC', \, \, n \in \bbN,
\\ %
& \ol{M (\gamma \omega)^{n + 1}}^n \simeq M (\gamma B^n \omega), \,
n \in \bbN,
\\ %
& \ol{N (C, \omega)} \simeq L (B C), \, C \in \calC',
\\ %
& \ol{M (\omega) N (\omega) M (\gamma \omega)^n}^{n + 1} \simeq L
(B^{n + 1} \omega), \, n \in \bbN_0,
\\ %
& \ol{M (\gamma C) M (\omega) N (\omega) M (\gamma \omega)^{n -
1}}^{n + 1} \simeq L (B^{n + 1} C), C \in \calC', \, n \in \bbN,
\\ %
& \ol{M (\omega) M (\gamma \omega)^n}^n \simeq M (B^n \omega), \, n
\in \bbN,
\\ %
& \ol{M (\gamma C) M (\omega) M (\gamma \omega)^n}^{n + 1} \simeq M
(B^{n + 1} C), \, C \in \calC', \, n \in \bbN_0,
\\ %
& \ol{N (\omega) M (\gamma \omega)^n}^n \simeq N (B^n \omega), \, n
\in \bbN,
\\ %
& \ol{M (\gamma C) M (\omega) M (\gamma \omega)^n}^{n + 1} \simeq N
(B^{n + 1} C), \, C \in \calC', \, n \in \bbN_0,
\\ %
& \ol{M (\gamma \omega)^n}^n (\lambda) \simeq R (B, \lambda, n), \,
n \in \bbN, \, \lambda \in k^*,
\\ %
& \ol{N (\omega) M (\gamma \omega)^n}^{n + 1} \simeq Q (B^{n + 1}),
\, n \in \bbN_0,
\\ %
& \ol{M (\omega) M (\gamma \omega)^n}^{n + 1} \simeq M (B^n \omega
\alpha), \, n \in \bbN_0,
\\ %
& \ol{M (\gamma \omega)^n M (x)}^n \simeq M (\gamma B^{n - 1} \omega
\alpha \xi), \, n \in \bbN,
\\ %
& \ol{M (\gamma \omega)^n}^n (\infty) \simeq M (\gamma B^{n - 1}
\omega \alpha), \, n \in \bbN,
\\ %
& \ol{M (\omega) M (\gamma \omega)^n M (x)}^{n + 1} \simeq M (B^n
\omega \alpha \xi), \, n \in \bbN_0,
\\ %
& \ol{M (\omega_{j_p}) N (\omega_{j_p}) M (\gamma \omega)^m}^{m}
\simeq N (\omega_{j_p}, B_{j_p} B^{m - 1} \omega), \, p \in [1, r],
\, m \in \bbN,
\\ %
& \ol{M (\gamma \omega)^n M (\omega_{j_p}) N (\omega_{j_p}) M
(\gamma \omega)^m}^{m + n} \simeq N (B_{j_p} B^{n - 1} \omega,
B_{j_p} B^{m -
1} \omega),\\ %
& \qquad \, p \in [1, r], \, n, m \in \bbN, \, n < m,
\\ %
& \ol{M (\gamma C) M (\gamma \omega)^n M (\omega_{j_p}) N
(\omega_{j_p}) M (\gamma \omega)^m}^{m + n + 1}
\\ %
& \qquad \simeq N (B_{j_p} B^n C, B_{j_p} B^{m - 1} \omega), \, p
\in [1, r], \, C \in \calC', \, n, m \in \bbN_0, \, n < m,
\\ %
& \ol{M (\omega_{j_p}) N (\omega_{j_p}) M (\gamma \omega)^m M
(\gamma C)}^{m + 1} \simeq N (\omega_{j_p}, B_{j_p} B^m C),
\\ %
& \qquad p \in [1, r], \, C \in \calC', \, m \in \bbN_0,
\\ %
& \ol{M (\gamma \omega)^n M (\omega_{j_p}) N (\omega_{j_p}) M
(\gamma \omega)^m M (\gamma C)}^{m + n + 1}
\\ %
& \qquad \simeq N (B_{j_p} B^{n - 1} \omega, B_{j_p} B^m C), \, p
\in [1, r], \, C \in \calC', \, n, m \in \bbN, \, n \leq m,
\\ %
& \ol{M (\gamma C') M (\gamma \omega)^n M (\omega_{j_p}) N
(\omega_{j_p}) M (\gamma \omega)^m M (\gamma C'')}^{m + n + 1}
\\ %
& \qquad \simeq N (B_{j_p} B^n C', B_{j_p} B^m C''),
\\ %
& \qquad p \in [1, r], \, C', C'' \in \calC', \, n, m \in \bbN_0, \,
(n, C') < (m, C''),
\\ %
& \ol{M (\gamma \omega)^n M (\omega_j)}^n \simeq M (B_j B^{n - 1}
\omega), \, j \in [j_0 + 1, j_{r + 1} - 1], \, n \in \bbN,
\\ %
& \ol{M (\gamma C) M (\gamma \omega)^n M (\omega_j)}^{n + 1} \simeq
M (B_j B^n C), \, j \in [j_0 + 1, j_{r + 1} - 1], \, n \in \bbN_0,
\\ %
& \ol{M (\gamma \omega)^n N (\omega_{j_p})}^n \simeq N (B_{j_p} B^{n
- 1} \omega), \, p \in [1, r], \, n \in \bbN,
\\ %
& \ol{M (\gamma C) M (\gamma \omega)^n N (\omega_{j_p})}^{n + 1}
\simeq N (B_{j_p} B^n C), \, p \in [1, r], \, C \in \calC', \, n \in
\bbN_0,
\\ %
& \ol{M (\gamma \omega)^m M (\gamma_p \omega_{j_p})}^m \simeq M
(\gamma_p B_{j_p} B^{m - 1} \omega), \, p \in [1, r], \, m \in \bbN,
\\ %
& \ol{M (\gamma \omega)^m N (C, \omega_{j_p})}^m \simeq N (C,
B_{j_p} B^{m - 1} \omega), \, C \in \calC_p', \, p \in [1, r], \, m
\in \bbN,
\\ %
& \ol{M (\gamma C) M (\gamma \omega)^m M (\gamma_p \omega_{j_p})}^{m
+ 1} \simeq M (\gamma_p B_{j_p} B^m C),
\\ %
& \qquad C \in \calC', \, p \in [1, r], \, m \in \bbN_0,
\\ %
& \ol{M (\gamma C') M (\gamma \omega)^m N (C'', \omega_{j_p})}^{m +
1} \simeq N (C'', B_{j_p} B^m C'),
\\ %
& \qquad C' \in \calC', \, C'' \in \calC_p', \, p \in [1, r], \, m
\in \bbN_0,
\\ %
& \ol{M (\gamma \omega)^m}^m (0) \simeq M (\xi \gamma B^{m - 1}
\omega), \, m \in \bbN,
\\ %
& \ol{M (\gamma C) M (\gamma \omega)^m}^{m + 1} \simeq M (\xi \gamma
B^m C), \, C \in \calC', \, m \in \bbN_0,
\\ %
& \ol{M (\omega_{j_p}) N (\omega_{j_p})}^1 \simeq N (\omega_{j_p},
\omega_{j_p} \alpha), \, p \in [1, r],
\\ %
& \ol{M (\gamma \omega)^m M (\omega_{j_p}) N (\omega_{j_p})}^{m + 1}
\simeq N (B_{j_p} B^{m - 1} \omega, \omega_{j_p} \alpha), \, p \in
[1, r], \, m \in \bbN,
\\ %
& \ol{M (\omega_{j_p}) N (\omega_{j_p}) M (\gamma \omega)^n}^{n + 1}
\simeq N (\omega_{j_p}, B_{j_p} B^{n - 1} \omega \alpha), \, p \in
[1, r], \, n \in \bbN,
\\ %
& \ol{M (\gamma \omega)^m M (\omega_{j_p}) N (\omega_{j_p}) M
(\gamma \omega)^n}^{n + m + 1} \simeq N (B_{j_p} B^{m - 1} \omega,
B_{j_p} B^{n - 1} \omega \alpha),
\\ %
& \qquad p \in [1, r], \, n, m \in \bbN,
\\ %
& \ol{M (\gamma C) M (\gamma \omega)^m M (\omega_{j_p}) N
(\omega_{j_p})}^{m + 2} \simeq N (B_{j_p} B^m C, \omega_{j_p}
\alpha),
\\ %
& \qquad C \in \calC', \, p \in [1, r], \, m \in \bbN_0,
\\ %
& \ol{M (\gamma C) M (\gamma \omega)^m M (\omega_{j_p}) N
(\omega_{j_p}) M (\gamma \omega)^n}^{n + m + 2}
\\ %
& \qquad \simeq N (B_{j_p} B^m C, B_{j_p} B^{n - 1} \omega \alpha),
\, C \in \calC', \, p \in [1, r], \, n \in \bbN, \, m \in \bbN_0,
\\ %
& \ol{M (\omega_{j_p}) N (\omega_{j_p}) M (x)}^1 \simeq N
(\omega_{j_p}, \omega_{j_p} \alpha \xi), \, p \in [1, r],
\\ %
& \ol{M (\gamma \omega)^m M (\omega_{j_p}) N (\omega_{j_p}) M
(x)}^{m + 1} \simeq N (B_{j_p} B^{m - 1} \omega, \omega_{j_p} \alpha
\xi),
\\ %
& \qquad p \in [1, r], \, m \in \bbN,
\\ %
& \ol{M (\omega_{j_p}) N (\omega_{j_p}) M (\gamma \omega)^n M
(x)}^{n + 1} \simeq N (\omega_{j_p}, B_{j_p} B^{n - 1} \omega \alpha
\xi),
\\ %
& \qquad p \in [1, r], \, n \in \bbN,
\\ %
& \ol{M (\gamma \omega)^m M (\omega_{j_p}) N (\omega_{j_p}) M
(\gamma \omega)^n M (x)}^{n + m + 1}
\\ %
& \qquad \simeq N (B_{j_p} B^{m - 1} \omega, B_{j_p} B^{n - 1}
\omega \alpha \xi), \, p \in [1, r], \, n, m \in \bbN,
\\ %
& \ol{M (\gamma C) M (\gamma \omega)^m M (\omega_{j_p}) N
(\omega_{j_p}) M (x)}^{m + 2} \simeq N (B_{j_p} B^m C, \omega_{j_p}
\alpha \xi),
\\ %
& \qquad C \in \calC', \, p \in [1, r], \, m \in \bbN_0,
\\ %
& \ol{M (\gamma C) M (\gamma \omega)^m M (\omega_{j_p}) N
(\omega_{j_p}) M (\gamma \omega)^n M (x)}^{n + m + 2}
\\ %
& \qquad \simeq N (B_{j_p} B^m C, B_{j_p} B^{n - 1} \omega \alpha
\xi), \, C \in \calC', \, p \in [1, r], \, n \in \bbN, \, m \in
\bbN_0,
\\ %
& \ol{M (\gamma_p \omega_{j_p})}^1 \simeq M (\gamma_p \omega_{j_p}
\alpha), \, p \in [1, r],
\\ %
& \ol{M (\gamma \omega)^n M (\gamma_p \omega_{j_p})}^{n + 1} \simeq
M (\gamma_p B_{j_p} B^{n - 1} \omega \alpha), \, p \in [1, r], \, n
\in \bbN,
\\ %
& \ol{M (\omega_j)}^1 \simeq M (\omega_j \alpha), \, j \in [j_0 + 1,
j_{r + 1} - 1],
\\ %
& \ol{M (\gamma \omega)^n M (\omega_j)}^{n + 1} \simeq M (B_j B^{n -
1} \omega \alpha), \, j \in [j_0 + 1, j_{r + 1} - 1], \, n \in \bbN,
\\ %
& \ol{N (C, \omega_{j_p})}^1 \simeq N (C, \omega_{j_p} \alpha), \, C
\in \calC_p, \, p \in [1, r],
\\ %
& \ol{M (\gamma \omega)^n N (C, \omega_{j_p})}^{n + 1} \simeq N (C,
B_{j_p} B^{n - 1} \omega \alpha), \, C \in \calC_p, \, p \in [1, r],
\, n \in \bbN,
\\ %
& \ol{M (\gamma_p \omega_{j_p}) M (x)}^1 \simeq M (\gamma_p
\omega_{j_p} \alpha \xi), \, p \in [1, r],
\\ %
& \ol{M (\gamma \omega)^n M (\gamma_p \omega_{j_p}) M (x)}^{n + 1}
\simeq M (\gamma_p B_{j_p} B^{n - 1} \omega \alpha \xi), \, p \in
[1, r], \, n \in \bbN,
\\ %
& \ol{M (\omega_j) M (x)}^1 \simeq M (\omega_j \alpha \xi), \, j \in
[j_0 + 1, j_{r + 1} - 1],
\\ %
& \ol{M (\gamma \omega)^n M (\omega_j) M (x)}^{n + 1} \simeq M (B_j
B^{n - 1} \omega \alpha \xi),
\\ %
& \qquad j \in [j_0 + 1, j_{r + 1} - 1], \, n \in \bbN,
\\ %
& \ol{N (C, \omega_{j_p}) M (x)}^1 \simeq N (C, \omega_{j_p} \alpha
\xi), \, C \in \calC_p, \, p \in [1, r],
\\ %
& \ol{M (\gamma \omega)^n N (C, \omega_{j_p}) M (x)}^{n + 1} \simeq
N (C, B_{j_p} B^{n - 1} \omega \alpha \xi),
\\ %
& \qquad C \in \calC_p, \, p \in [1, r], \, n \in \bbN,
\\ %
& \ol{0}^1 \simeq M (z),
\\ %
& \ol{M (\gamma \omega)^n}^{n + 1} \simeq M (\xi \gamma B^{n - 1}
\omega \alpha), \, n \in \bbN,
\\ %
& \ol{M (x)}^1 \simeq M (\xi),
\\ %
& \ol{M (\gamma \omega)^n M (x)}^{n + 1} \simeq M (\xi \gamma B^{n
-1} \omega \alpha \xi), \, n \in \bbN,
\\ %
& \ol{M (\gamma \omega)^n M (\omega_{j_p}) N (\omega_{j_p})}^{n + 2}
\simeq N (B_{j_p} B^{n - 1} \omega \alpha, \omega_{j_p} \alpha), \,
p \in [1, r], \, n \in \bbN,
\\ %
& \ol{M (\gamma \omega)^n M (\omega_{j_p}) N (\omega_{j_p}) M
(\gamma \omega)^m}^{n + m + 2} \simeq N (B_{j_p} B^{n - 1} \omega
\alpha, B_{j_p} B^{m - 1} \omega \alpha),
\\ %
& \qquad p \in [1, r], \, n, m \in \bbN, \, m < n,
\\ %
& \ol{M (x) M (\omega_{j_p}) N (\omega_{j_p})}^{2} \simeq N
(\omega_{j_p} \alpha \xi, \omega_{j_p} \alpha), \, p \in [1, r],
\\ %
& \ol{M (x) M (\gamma \omega)^n M (\omega_{j_p}) N
(\omega_{j_p})}^{n + 2} \simeq N (B_{j_p} B^{n - 1} \omega \alpha
\xi, \omega_{j_p} \alpha),
\\ %
& \qquad p \in [1, r], \, n \in \bbN,
\\ %
& \ol{M (x) M (\gamma \omega)^n M (\omega_{j_p}) N (\omega_{j_p}) M
(\gamma \omega)^m}^{n + m + 2}
\\ %
& \qquad \simeq N (B_{j_p} B^{n - 1} \omega \alpha \xi, B_{j_p} B^{m
- 1} \omega \alpha), \, p \in [1, r], \, n, m \in \bbN, \, m \leq n,
\\ %
& \ol{M (\gamma \omega)^n M (\omega_{j_p}) N (\omega_{j_p}) M
(x)}^{n + 2} \simeq N (B_{j_p} B^{n - 1} \omega \alpha, \omega_{j_p}
\alpha \xi),
\\ %
& \qquad p \in [1, r], \, n \in \bbN,
\\ %
& \ol{M (\gamma \omega)^n M (\omega_{j_p}) N (\omega_{j_p}) M
(\gamma \omega)^m M (x)}^{n + m + 2}
\\ %
& \qquad \simeq N (B_{j_p} B^{n - 1} \omega \alpha, B_{j_p} B^{m -
1} \omega \alpha \xi), \, p \in [1, r], \, n, m \in \bbN, \, m < n,
\\ %
& \ol{M (x) M (\gamma \omega)^n M (\omega_{j_p}) N (\omega_{j_p}) M
(x)}^{n + 2} \simeq N (B_{j_p} B^{n - 1} \omega \alpha \xi,
\omega_{j_p} \alpha \xi),
\\ %
& \qquad p \in [1, r], \, n \in \bbN,
\\ %
& \ol{M (x) M (\gamma \omega)^n M (\omega_{j_p}) N (\omega_{j_p}) M
(\gamma \omega)^m M (x)}^{n + m + 2}
\\ %
& \qquad \simeq N (B_{j_p} B^{n - 1} \omega \alpha \xi, B_{j_p} B^{m
- 1} \omega \alpha \xi), \, p \in [1, r], \, n, m \in \bbN, \, m <
n,
\\ %
& \ol{N (\omega_{j_p})}^1 \simeq N (\omega_{j_p} \alpha), \, p \in
[1, r],
\\ %
& \ol{N (\gamma \omega)^n M (\omega_{j_p})}^1 \simeq N (B_{j_p} B^{n
- 1} \omega \alpha), \, p \in [1, r], \, n \in \bbN,
\\ %
& \ol{N (\omega_{j_p}) M (x)}^1 \simeq N (\omega_{j_p} \alpha \xi),
\, p \in [1, r],
\\ %
& \ol{N (\gamma \omega)^n M (\omega_{j_p}) M (x)}^1 \simeq N
(B_{j_p} B^{n - 1} \omega \alpha \xi), \, p \in [1, r], \, n \in
\bbN,
\end{align*}

\subsection{}
We now prove Lemma~\ref{lemmR}. Let $x'$ be an admissible index of
$Q^x$. Let first $x' = x_{i, j}$, $x' \neq z$. In this case
$R_{x'}^x = R_{x'}$ . If $i \neq i_0$ or $i = i_0$ and $j < j_0$
then also $\Hom (R_{x'}^x, \mod A^x) = \Hom (R_{x'}, \mod A)$. If $i
= i_0$ and $j_0 < j$, then the indecomposable objects of $\Hom
(R_{x'}^x, \mod A^x)$ are the indecomposable objects of $\Hom
(R_{x'}, \mod A)$ and
\begin{gather*}
\ol{M (\gamma \omega)^n M (\omega_{j - 1})}^n, \, \ol{M (\gamma
\omega)^n M (\omega_j)}^n, \, n \in \bbN,
\\ %
\ol{M (\gamma C) M (\gamma \omega)^n M (\omega_{j - 1})}^{n + 1},
\, \ol{M (\gamma C) M (\gamma \omega)^n M (\omega_j)}^{n + 1}, C
\in \calC', n \in \bbN_0,
\\ %
\ol{M (\gamma \omega)^n M (\omega_{j - 1})}^{n + 1}, \, \ol{M
(\gamma \omega)^n M (\omega_j)}^{n + 1}, \, n \in \bbN_0,
\\ %
\intertext{and}
\ol{M (\gamma \omega)^n M (\omega_{j - 1}) M (x)}^{n + 1}, \,
\ol{M (\gamma \omega)^n M (\omega_j) M (x)}^{n + 1}, \, n \in
\bbN_0.
\end{gather*}

Assume now that $x' = z_{i, j}$. If $i \neq i_0$, then again
$R_{x'}^x = R_{x'}$ and $\Hom (R_{x'}^x, \mod A^x) = \Hom (R_{x'},
\mod A)$. If $i = i_0$ then $j = j_p$ for $p \in [1, r]$. Moreover,
$R_{x'}^x = \ol{R_{x'}}$ and the indecomposable objects of $\Hom
(R_{x'}^x, \mod A^x)$ are
\begin{gather*}
M (\gamma_p C), \, C \in \calC_p, \, M (x'),
\\ %
\ol{N (C, \omega_{j_q})}, \, \ol{M (\gamma_q \omega_{j_q})}, \,
\ol{N (\omega_{j_q})}, C \in \calC_{j_q}', \, \, q \in [p, r],
\\ %
\ol{M (\omega_l)}, \, l \in [j_p, \ldots, j_{r + 1} - 1], \,
\ol{0},
\\ %
\ol{M (\gamma \omega)^n M (\gamma_p \omega_{j_p})}^n, \, \ol{M
(\gamma \omega)^{n - 1} M (\omega_{j_q}) N (\omega_{j_q})}^n, \, n
\in \bbN, \, q \in [p, r],
\\ %
\ol{M (\gamma C) M (\gamma \omega)^n M (\gamma_p \omega_{j_p})}^{n
+ 1}, \, C \in \calC', \, n \in \bbN_0,
\\ %
\ol{M (\gamma C) M (\gamma \omega)^n M (\omega_{j_q}) N
(\omega_{j_q})}^{n + 2}, \, C \in \calC_{j_q}', \, n \in \bbN_0,
\, q \in [p, r],
\\ %
\ol{M (\gamma \omega)^n M (\gamma_p \omega_{j_p})}^{n + 1}, \,
\ol{M (\gamma \omega)^{n - 1} M (\gamma_p \omega_{j_p}) M (x)}^n,
\, n \in \bbN,
\\ %
\ol{M (\gamma \omega)^n M (\omega_{j_q}) N (\omega_{j_q})}^{n +
2}, \, n \in \bbN_0, \, q \in [p, r],
\\ %
\intertext{and}
\ol{M (\gamma \omega)^n M (\omega_{j_q}) N (\omega_{j_q}) M
(x)}^{n + 2}, \, n \in \bbN_0, \, q \in [p, r].
\end{gather*}

Finally, assume that $x' = z$ (it is possible, if $p_{i_0} +
|T_{i_0}| \not \in S_{i_0}$). In this case $R_{x'}^x = \ol{X_x M
(x_{i_0, p_{i_0} + |T_{i_0}|})}$. It follows that the indecomposable
objects of $\Hom (R_{x'}^x, \mod A^x)$ are
\begin{gather*}
\ol{M (C) M (x_{i_0, p_{i_0} + |T_{i_0}|})}, \, \ol{M (C)}, \, C
\in \calC_x,
\\ %
\intertext{and} \ol{M (x_{i_0, p_{i_0} + |T_{i_0}|})}, \, \ol{0}.
\end{gather*}

\subsection{}
Now we indicate how to prove Lemma~\ref{lemmX}. Let $x'$ be an
admissible index of $Q^x$. If $x' = x_{i, j}$, $x' \neq z$, then
$X_{x'}^x = X_{x'}$. If in addition, $i \neq i_0$ or $i = i_0$ and
$j < j_0$, then $\Hom (X_{x'}^x, \mod A^x) = \Hom (X_{x'}, \mod A)$.
Let $i = i_0$ and $j_0 < j$. Then the indecomposable objects of
$\Hom (X_{x'}^x, \mod A^x)$ are the indecomposable objects of $\Hom
(X_{x'}, \mod A)$ and
\begin{gather*}
\ol{M (\gamma \omega)^{n + 1} M (\omega_j)}^{n + 1}, \, \ol{M
(\gamma C) M (\gamma \omega)^n M (\omega_j)}^{n + 1}, C \in
\calC', n \in \bbN_0,
\\ %
\intertext{and}
\ol{M (\gamma \omega)^n M (\omega_j)}^{n + 1}, \, \ol{M (\gamma
\omega)^n M (\omega_j) M (x)}^{n + 1}, \, n \in \bbN_0.
\end{gather*}

Assume now that $x' = z_{i, j}$. Again $X_{x'}^x = X_{x'}$ and if $i
\neq i_0$ then $\Hom (X_{x'}^x, \mod A^x) = \Hom (X_{x'}, \mod A)$.
Let $i = i_0$. Then $j = j_p$ for $p \in [1, r]$. The indecomposable
objects of $\Hom (X_{x'}^x, \mod A^x)$ are the indecomposable
objects of $\Hom (X_{x'}, \mod A)$ and
\begin{gather*}
\ol{M (\gamma \omega)^{n + 1} M (\gamma_p \omega_{j_p})}^{n + 1},
\, \ol{M (\gamma C) M (\gamma \omega)^n M (\gamma_p
\omega_{j_p})}^{n + 1}, \, C \in \calC', \, n \in \bbN_0,
\\ %
\\ %
\intertext{and}
\ol{M (\gamma \omega)^n M (\gamma_p \omega_{j_p})}^{n + 1}, \,
\ol{M (\gamma \omega)^n M (\gamma_p \omega_{j_p}) M (x)}^{n + 1},
\, n \in \bbN_0.
\end{gather*}

Finally, let $x' = z$. In this case $X_{x'}^x = \ol{X_x}$ and the
indecomposable objects of $\Hom (X_{x'}, \mod A^x)$ are
\[
\ol{M (C)}, \, C \in \calC_x, \text{ and } \ol{0}.
\]

\subsection{}
It remains to give the proof of Lemma~\ref{lemmI}. Let $x' = x_{i,
j}$ be the vertex of $Q^x$ satisfying the hypothesis of
Lemma~\ref{lemmI}. If $i \neq i_0$ then $I_{x'}^x = I_{x'}$ and
$\Hom (I_{x'}^x, \mod A^x) = \Hom (I_{x'}, \mod A)$. Assume now that
$i = i_0$. Then $j \in [j_0 + 1, j_{r + 1}$. Let $p = \min \{ q \in
[1, r + 1] \mid j \leq j_q \}$. First consider the case $j \neq j_{r
+ 1}$. Then $I_{x'}^x = \ol{I_{x'}}$ and the indecomposable objects
of $\Hom (I_{x'}^x, \mod A^x)$ are
\begin{gather*}
\ol{N (C, \omega_{j_q})}, \, \ol{M (\gamma_q \omega_{j_q})}, \,
\ol{N (\omega_{j_q})}, \, C \in \calC_{j_q}', \, q \in [p, r], \,
\\ %
\ol{M (\omega_l)}, \, l \in [j, \ldots, j_{r + 1} - 1], \, \ol{0},
\\ %
\ol{M (\gamma \omega)^{n - 1} M (\omega_{j_q}) N
(\omega_{j_q})}^n, \, n \in \bbN, \, q \in [p, r],
\\ %
\ol{M (\gamma C) M (\gamma \omega)^n M (\omega_{j_q}) N
(\omega_{j_q})}^{n + 2}, \, C \in \calC_{j_q}', \, n \in \bbN_0,
\, q \in [p, r],
\\ %
\ol{M (\gamma \omega)^n M (\omega_{j_q}) N (\omega_{j_q})}^{n +
2}, \, n \in \bbN_0, \, q \in [p, r],
\\ %
\intertext{and}
\ol{M (x) M (\gamma \omega)^n M (\omega_{j_q}) N
(\omega_{j_q})}^{n + 2}, \, n \in \bbN_0, \, q \in [p, r].
\end{gather*}
If $j = j_{r + 1}$ then the claim is clear.



\end{document}